\documentclass[final]{siamltex}
\usepackage[vlined,algosection,ruled]{algorithm2e}
\usepackage{amsfonts}
\usepackage{amsmath,amssymb}
\usepackage{comment}
\usepackage{tikz}
\usepackage{paralist}
\usepackage{pgfplots}
\usepackage{url}

\def\gmresk{{\rm GMRES($k$)}} 

\title{A parallel sweeping preconditioner for\\ 
       heterogeneous 3D Helmholtz equations\thanks{This work was partially 
       supported by the sponsors of the Texas Consortium for Computational 
       Seismology.}}

\author{Jack Poulson\thanks{ICES,
        University of Texas at Austin, 1 University Station C0200, 
        Austin, TX, 78712 (jack.poulson@gmail.com). This author was also 
        supported by a CAM fellowship.} \and 
        Bj\"orn Engquist\thanks{Department of Mathematics and ICES,
        University of Texas at Austin, 1 University Station C1200,
        Austin, TX, 78712 (engquist@ices.utexas.edu). This author was also 
        supported by NSF grant DMS-1016577.} \and 
        Siwei Li\thanks{Jackson School of Geosciences,
        University of Texas at Austin, 1 University Station C1160,
        Austin, TX, 78712 (siwei.li@utexas.edu).} \and
        Lexing Ying\thanks{Department of Mathematics and ICES, 
        University of Texas at Austin, 1 University Station C1200,
        Austin, TX, 78712 (lexing@math.utexas.edu). This author was supported 
        by NSF CAREER grant DMS-0846501, NSF grant DMS-1016577, and 
        funding from KAUST.}}

\begin{document}

\maketitle

\begin{abstract}
A parallelization of a sweeping preconditioner for 3D Helmholtz equations 
without large cavities is introduced and benchmarked for several 
challenging velocity models. 
The setup and application costs of the
sequential preconditioner are shown to be $O(\gamma^2 N^{4/3})$ 
and $O(\gamma N \log N)$, where $\gamma(\omega)$ denotes the modestly 
frequency-dependent number of grid points per Perfectly Matched Layer. 
Several computational and memory improvements are introduced relative to 
using black-box sparse-direct solvers for the auxiliary problems, and 
competitive runtimes and iteration counts are reported for high-frequency 
problems distributed over thousands of cores.
Two open-source packages are released along with this paper: {\em Parallel 
Sweeping Preconditioner (PSP)} and the underlying distributed multifrontal 
solver, {\em Clique}.
\end{abstract}



\pagestyle{myheadings}
\thispagestyle{plain}
\markboth{J. POULSON ET AL.}{PARALLEL SWEEPING PRECONDITIONER}

\section{Introduction}
\label{section:intro}
While definite elliptic partial differential equations can be efficiently 
solved by a wide variety of techniques (e.g., multigrid, ILU, or structured 
matrix factorizations), indefinite elliptic equations tend to be more 
challenging. 
This paper is concerned with three-dimensional heterogeneous Helmholtz 
equations of the form,
\begin{equation}\label{helmholtz}
  \mathcal{A} u \equiv 
  \left[-\Delta - \frac{\omega^2}{c^2(x)}\right] u(x) =f(x),
\end{equation}
where $c(x)$ is the spatially varying wave speed, and $u(x) e^{-i \omega t}$ 
is the time-harmonic response to an acoustic wave equation with forcing 
function $f(x) e^{-i \omega t}$.
It is important to recognize that $-\Delta$ is positive-definite and that 
its combination with the negative-definite $-\frac{\omega^2}{c^2}$ term 
results in an indefinite system.

Before discussing the overall asymptotic complexity of solution techniques,
it is helpful to first motivate why high frequency problems require large 
numbers of degrees of freedom: Given the wave speed bounds 
$c_{\text{min}} \le c(x) \le c_{\text{max}}$, we can define the minimum 
wavelength as $\lambda_{\text{min}} = 2 \pi c_{\text{min}}/\omega$.
In order to resolve oscillations in the solution using piecewise polynomial 
basis functions, e.g., with finite-difference and finite-element methods, it is 
necessary to increase the number of degrees of freedom in each direction at 
least linearly with the number of wavelengths spanned by the domain. In order
to combat pollution effects~\cite{Babuska-pollution}, which are closely related 
to phase errors in the discrete solution, one must use asymptotically more
than a constant number of grid points per wavelength with standard 
discretization schemes. Nevertheless, it is common practice to resolve the 
domain to as few as ten points per wavelength. 
In any case, piecewise polynomial discretizations require $\Omega(\omega^d)$
degrees of freedom in $d$ dimensions.

Until recently, doubling the frequency of Eq.\ \eqref{helmholtz} not only 
increased the size of the linear system by at least a factor of $2^d$, it also 
doubled the number of iterations required for convergence with standard 
preconditioned Krylov methods~\cite{BollhoeferGroteSchenk-ILU,Erlangga-advances,ErnstGander-classical}. 
Thus, denoting the number of degrees of freedom in a three-dimensional 
finite-element or finite-difference discretization as $N = \Omega(\omega^3)$,
every linear solve required $\Omega(\omega^4)$ work with traditional iterative 
techniques.
Engquist and Ying recently introduced two classes of {\em sweeping} 
preconditioners for Helmholtz equations without large 
cavities~\cite{EngquistYing-HMat,EngquistYing-PML}: Both approaches 
approximate a block $LDL^T$ factorization of the Helmholtz operator in 
block tridiagonal form in a manner which exploits a radiation boundary 
condition. 
The first approach performs a block tridiagonal factorization algorithm in 
$\mathcal{H}$-matrix arithmetic~\cite{Hackbusch-HMat,Grasedyck-HMat}, 
while the second approach approximates the Schur complements of the 
factorization using auxiliary problems with artificial radiation boundary 
conditions. 
Though the $\mathcal{H}$-matrix sweeping preconditioner has theoretical 
support for two-dimensional 
problems~\cite{EngquistYing-HMat,MartinssonRokhlin-elongated}, 
there is not yet justification for three-dimensional problems. 

This paper therefore focuses on the second approach, which relies on 
multifrontal factorizations~\cite{Liu-review,Schreiber-etree,DuffReid-multifrontal,George-nested} 
of the approximate auxiliary problems in order to achieve an 
$O(\gamma^2 N^{4/3})$ setup cost and an $O(\gamma N \log N)$ application 
cost, where $\gamma(\omega)$ denotes the number of grid points used for each 
Perfectly Matched Layer (PML)~\cite{Johnson-PML}. While the sweeping 
preconditioner is competitive with existing techniques even for a single 
right-hand side, its main advantage is for problems with large numbers of 
right-hand sides, as the preconditioner appears to converge in $O(1)$ 
iterations for problems without large cavities. Thus, after setting up 
the preconditioner, typically only $O(\gamma N \log N)$ work is required for 
each solution.

\subsection{{\em Moving PML} sweeping preconditioner}
The focus of this paper is on parallelization of the second class of sweeping
preconditioners mentioned above, which makes use of auxiliary problems with 
artificial radiation boundary conditions in order to approximate the Schur 
complements that arise during block $LDL^T$ factorization. The approach is 
referred to as a {\em moving PML} preconditioner since the introductory paper
represented the artificial radiation boundary conditions using PML.

One interpretation of radiation conditions is that they allow for the 
analysis of a finite portion of an infinite domain, as their purpose is to 
enforce the condition that waves propagate outwards and not reflect at the 
boundary of the truncated domain.
This concept is crucial to understanding the Schur complement approximations 
that take place within the moving PML sweeping preconditioner which is 
reintroduced in this paper for the sake of completeness.

For the sake of simplicity, we will describe the preconditioner in the context
of a second-order finite-difference discretization over the unit cube,
with PML used to approximate a radiation boundary condition on the $x_3=0$ face
and homogeneous Dirichlet boundary conditions applied on all other boundaries
(an $x_1 x_3$ cross-section is shown in Fig.~\ref{fig:plane-with-single-pml}). 
If the domain is discretized into an 
$n \times n \times n$ grid, then ordering the vertices in the grid such that 
vertex $(i_1,i_2,i_3)$ is assigned index 
$i_1 + i_2 n + i_3 n^2$ results in a block tridiagonal system of equations, say
\begin{equation}
  \left(\begin{array}{ccccc} 
          A_{0,0} & A_{1,0}^T &        &           &               \\
          A_{1,0} & A_{1,1}   & \ddots &           &               \\
                  & \ddots    & \ddots & \ddots    &               \\
                  &           & \ddots & \ddots    & A_{n-1,n-2}^T \\
                  &           &        & A_{n-1,n-2} & A_{n-1,n-1}
  \end{array}\right) 
  \left(\begin{array}{c}u_0\\ u_1\\\vdots\\u_{n-2}\\u_{n-1}\end{array}\right)=
  \left(\begin{array}{c}f_0\\ f_1\\\vdots\\f_{n-2}\\f_{n-1}\end{array}\right),
\end{equation}
where $A_{i,j}$ propagates sources from the $i$'th $x_1x_2$ plane into the 
$j$'th $x_1x_2$ plane, and the overall linear system is complex symmetric 
({\em not} Hermitian) due to the imaginary terms introduced by the 
PML~\cite{EngquistYing-PML}. 

\begin{figure}
\centering
\begin{tikzpicture}[scale=4.1]
\tikzstyle{every node}=[font=\footnotesize]

\draw[->,thick] (-0.05,0) -- (-0.05,0.2);
\draw (-0.05,0.25) node {$x_1$};

\draw[->,thick] (0,-0.05) -- (0.2,-0.05);
\draw (0.25,-0.05) node {$x_3$};

\shade[left color=black!87!,inner color=gray,right color=white]
  (0,0) rectangle (0.125,1);

\draw[step=0.03125] (0,0) grid (1,1);

\node[fill=gray!5!,text width=1.3cm] (note1) at (0.5625,0.5)
{region of interest};

\draw (0.0625,0.5   ) node[fill=gray!5!] {PML};

\end{tikzpicture}
\caption{An $x_1 x_3$ cross section of a cube with PML on its $x_3=0$ face.
The domain is shaded in a manner which loosely corresponds to its extension into
the complex plane.}
\label{fig:plane-with-single-pml}
\end{figure}

If we were to ignore the sparsity within each block, then the following 
na\"ive factorization and solve algorithms would be appropriate for a 
direct solver, where the application of $S_i^{-1}$ implicitly makes use of the 
factorization of $S_i$.

\begin{algorithm}
\DontPrintSemicolon
$S_0 := A_{0,0}$\;
Factor $S_0$\;
\For{$i=0,...,n-2$}{
  $S_{i+1} := A_{i+1,i+1} - A_{i+1,i} S_i^{-1} A_{i+1,i}^T$\;
  Factor $S_{i+1}$\;
}
\caption{Na\"ive block tridiagonal $LDL^T$ factorization. 
$O(n^7)=O(N^{7/3})$ work is required.}
\label{alg:naive-fact}
\end{algorithm}

\begin{algorithm}
\DontPrintSemicolon
\tcp{Apply $L^{-1}$}
\For{$i=0,...,n-2$}{
  $u_{i+1} := u_{i+1} - A_{i+1,i}(S_i^{-1} u_i)$\;
}
\tcp{Apply $D^{-1}$}
\For{$i=0,...,n-1$}{
  $u_i := S_i^{-1} u_i$\;
}
\tcp{Apply $L^{-T}$}
\For{$i=n-2,...,0$}{
  $u_i := u_i - S_i^{-1}(A_{i+1,i}^T u_{i+1})$\;
}
\caption{Na\"ive block $LDL^T$ solve. $O(n^5)=O(N^{5/3})$ work is required.}
\label{alg:naive-solve}
\end{algorithm}
While the computational complexities of 
Algs.~\ref{alg:naive-fact} and \ref{alg:naive-solve} are significantly higher 
than multifrontal algorithms
\cite{Liu-review,DuffReid-multifrontal,George-nested},
which have $O(N^2)$ factorization complexity and $O(N^{4/3})$ 
solve complexity for regular three-dimensional meshes, they are the conceptual 
starting points of the sweeping preconditioner.\footnote{In fact, they are the 
starting points of {\em both} classes of sweeping preconditioners. 
The $\mathcal{H}$-matrix approach essentially executes these algorithms with 
$\mathcal{H}$-matrix arithmetic.} 

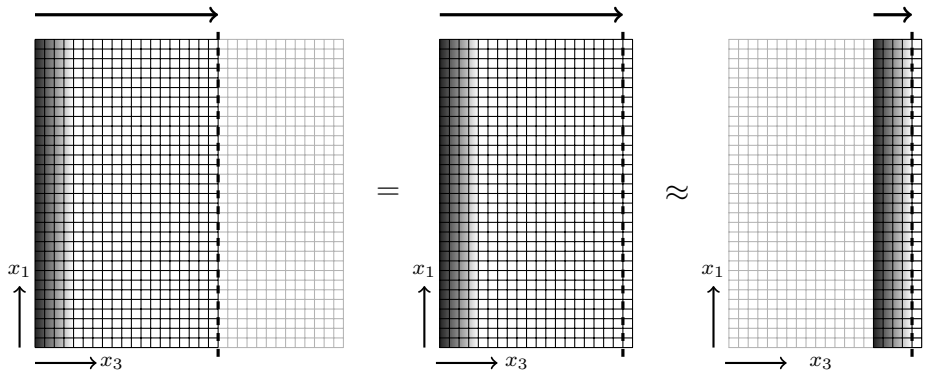
\begin{figure}
\centering
\begin{tikzpicture}[scale=4.1]
\tikzstyle{every node}=[font=\footnotesize]

%
%

\draw[->,very thick] (0,1.08) -- (0.59375,1.08);

\draw[->,thick] (-0.05,0) -- (-0.05,0.2);
\draw (-0.05,0.25) node {$x_1$};
\draw[->,thick] (0,-0.05) -- (0.2,-0.05);
\draw (0.25,-0.05) node {$x_3$};

\shade[left color=black!87!,inner color=gray,right color=white] 
  (0,0) rectangle (0.125,1);

\draw[step=0.03125] (0,0) grid (0.59375,1);

\draw[step=0.03125,opacity=0.3] (0.59375,0) grid (1,1);

\draw[very thick,dashed] (0.59375,-0.03) -- (0.59375,1.03);


%
%

\draw (1.145,0.5) node {\large $=$};

%
%

\draw[->,very thick] (1.3125,1.08) -- (1.90625,1.08);

\draw[->,thick] (1.2625,0) -- (1.2625,0.2);
\draw (1.2625,0.25) node {$x_1$};
\draw[->,thick] (1.3,-0.05) -- (1.5,-0.05);
\draw (1.5625,-0.05) node {$x_3$};

\shade[left color=black!87!,inner color=gray,right color=white] 
  (1.3125,0) rectangle (1.4375,1);

\draw[step=0.03125] (1.3120,0) grid (1.9375,1);

\draw[very thick,dashed] (1.90625,-0.03) -- (1.90625,1.03);


%
%

\draw (2.0825,0.5) node {\large $\approx$};

%
%

\draw[->,very thick] (2.71875,1.08) -- (2.84375,1.08);

\draw[->,thick] (2.2,0) -- (2.2,0.2);
\draw (2.2,0.25) node {$x_1$};
\draw[->,thick] (2.2375,-0.05) -- (2.4375,-0.05);
\draw (2.55,-0.05) node {$x_3$};

\draw[opacity=0.3,step=0.03125] (2.249,0) grid (2.8755,1);

\shade[left color=black!87!,inner color=gray,right color=white] 
  (2.71875,0) rectangle (2.84375,1);

\draw[step=0.03125] (2.7187,0) grid (2.875,1);

\draw[very thick,dashed] (2.84375,-0.03) -- (2.84375,1.03);


\end{tikzpicture}
\caption{(Left) A depiction of the portion of the domain involved in the 
computation of the Schur complement of an $x_1 x_2$ plane 
(marked with the dashed line)
with respect to all of the planes to its left during execution of 
Alg.~\ref{alg:naive-fact}. (Middle) An equivalent auxiliary
problem which generates the same Schur complement; the original domain is 
truncated just to the right of the marked plane and a homogeneous Dirichlet 
boundary condition is placed on the cut.
(Right) A local auxiliary problem for generating an approximation to the 
relevant Schur complement; the radiation boundary condition of the exact
auxiliary problem is moved next to the marked plane.}
\label{fig:auxiliary}
\end{figure}

Suppose that we paused Alg.~\ref{alg:naive-fact} after computing the $i$'th 
Schur complement, $S_i$, where the $i$'th $x_1 x_2$ plane is in the
interior of the domain (i.e., it is not in a PML region).
Due to the ordering imposed on the degrees of freedom of the discretization, 
the first $i$ Schur complements are equivalent to those that would have been 
produced from applying Alg.~\ref{alg:naive-fact} to an auxiliary problem 
formed by placing a homogeneous Dirichlet boundary condition on the 
$(i+1)$'th $x_1x_2$ plane and ignoring all of the 
successive planes (see the left half of Fig.~\ref{fig:auxiliary}). 
While this observation does not immediately yield any computational savings, 
it does allow for a qualitative description of the inverse of the $i$'th Schur 
complement, $S_i^{-1}$: it is the restriction of the half-space Green's function
of the exact auxiliary problem onto the $i$'th $x_1x_2$ plane 
(recall that PML can be interpreted as approximating the effect of a domain 
extending to infinity).

The main approximation made in the sweeping preconditioner can now be 
succinctly described: since $S_i^{-1}$ is a restricted half-space Green's 
function which incorporates the velocity field of the first $i$ planes, it is 
natural to approximate it with another restricted half-space
Green's function which only takes into account the {\em local} velocity field,
i.e., by moving the PML next to the $i$'th plane (see the right half of 
Fig.~\ref{fig:auxiliary}). 

If we use $\gamma(\omega)$ to denote the number of grid points of PML as a 
function of the frequency, $\omega$, then it is 
important to recognize that the depth of the approximate auxiliary problems
in the $x_3$ direction is $\Omega(\gamma)$.\footnote{In all of the experiments 
in this paper, $\gamma(\omega)$ was either 5 or 6, and the subdomain depth 
never exceeded 10.}
If the depth of the approximate auxiliary problems was $O(1)$, then 
combining nested dissection with the multifrontal method over the regular 
$n \times n \times O(1)$ mesh would only require $O(n^3)$
work and $O(n^2 \log n)$ storage \cite{George-nested}. 
However, if the PML size is a slowly growing function of frequency, then 
applying 2D nested dissection to the {\em quasi-2D} domain requires 
$O(\gamma^3 n^3)$ work and $O(\gamma^2 n^2 \log n)$ storage,
as the number of degrees of freedom in each front increases by a factor of 
$\gamma$ and dense factorizations have cubic complexity.

Let us denote the quasi-2D discretization of the local auxiliary problem for
the $i$'th plane as $H_i$, and its corresponding approximation to the Schur 
complement $S_i$ as $\tilde S_i$. Since $\tilde S_i$ is essentially a dense 
matrix, it is advantageous to come up with 
an abstract scheme for applying $\tilde S_i^{-1}$: Assuming that $H_i$ was 
ordered in a manner consistent with the 
$(i_1,i_2,i_3) \mapsto i_1+i_2 n+i_3 n^2$ global ordering, the degrees of 
freedom corresponding to the $i$'th plane come last and we find that
\begin{equation}
  H_i^{-1}=\left(\begin{array}{cc} \star & \star \\ \star & \tilde S_i^{-1}
           \end{array}\right),
\end{equation}
where the irrelevant portions of the inverse have been marked with a $\star$.
Then, trivially,
\begin{equation}
  H_i^{-1} \left(\begin{array}{c}0\\ u_i\end{array}\right)= 
  \left(\begin{array}{cc}\star &\star\\ \star &\tilde S_i^{-1}\end{array}\right)
  \left(\begin{array}{c}0\\ u_i\end{array}\right)= 
  \left(\begin{array}{c}\star\\ \tilde S_i^{-1}u_i\end{array}\right),
\end{equation}
which implies a method for quickly computing $\tilde S_i^{-1} u_i$ given 
a factorization of $H_i$:

\begin{algorithm}
\DontPrintSemicolon
Form $\hat u_i$ as the extension of $u_i$ by zero over the artificial PML\;
Form $\hat v_i := H_i^{-1} \hat u_i$\;
Extract $\tilde S_i^{-1} u_i$ from the relevant entries of $\hat v_i$\;
\caption{Application of $\tilde S_i^{-1}$ to $u_i$ given a multifrontal
factorization of $H_i$. $O(\gamma^2 n^2 \log n)$ work is required.}
\label{alg:apply-inverse}
\end{algorithm}
{\em From now on, we write $T_i$ to refer to the application of the 
(approximate) inverse of the Schur complement for the $i$'th plane.} 

There are two more points to discuss before defining the full serial algorithm. 
The first is that, rather than performing an approximate $LDL^T$ factorization 
using a discretization of $\mathcal{A}$, the preconditioner is instead built 
from a discretization of an {\em artificially damped} version of the Helmholtz 
operator, say
\begin{equation}\label{artificial-helmholtz}
  \mathcal{J} \equiv 
  \left[-\Delta - \frac{(\omega+i\alpha)^2}{c^2(x)}\right],
\end{equation}
where $\alpha \approx 2 \pi$ is responsible for the artificial damping. This 
is in contrast to shifted Laplacian preconditioners 
\cite{BaylissGoldsteinTurkel-shifted,ErlanggaVuikOosterlee-class}, where 
$\alpha$ is typically $O(\omega)$ \cite{ErnstGander-classical}, and our 
motivation is to avoid introducing large long-range dispersion error by 
damping the long range interactions in the preconditioner. Just as $A$ 
refers to the discretization of the original Helmholtz operator, $\mathcal{A}$,
we will use $J$ to refer to the discretization of the artificially damped 
operator,
$\mathcal{J}$.

The second point is much easier to motivate: since the artificial PML in 
each approximate auxiliary problem is of depth $\gamma(\omega)$, processing 
a single plane at a time would require processing $O(n)$ subdomains with 
$O(\gamma^3 n^3)$ work each. Clearly, processing $O(\gamma)$ planes at once
has the same asymptotic complexity as processing a single plane, and so 
combining $O(\gamma)$ planes reduces the setup cost from $O(\gamma^3 N^{4/3})$ 
to $O(\gamma^2 N^{4/3})$.
Similarly, the memory usage becomes $O(\gamma N \log N)$ instead of 
$O(\gamma^2 N \log N)$.\footnote{Note that increasing the number of planes 
per panel provides a mechanism for interpolating between the sweeping 
preconditioner and a full multifrontal factorization.} If we refer to these 
sets of $O(\gamma)$ contiguous planes as {\em panels}, which we label from 
$0$ to $m-1$, where $m=O(n/\gamma)$, and correspondingly redefine 
$\{u_i\}$, $\{f_i\}$, $\{S_i\}$, $\{T_i\}$, and $\{H_i\}$, 
we have the following algorithm for setting up an approximate 
block $LDL^T$ factorization of the discrete artificially damped Helmholtz 
operator:

\begin{algorithm}
\DontPrintSemicolon
$S_0 := J_{0,0}$\;
Factor $S_0$\;
\For{$i=1,...,m-1$}{
  Form $H_i$ by prefixing PML to $J_{i,i}$\;
  Factor $H_i$\;
}
\caption{Setup phase of the sweeping preconditioner. 
$O(\gamma^2 N^{4/3})$ work is required.}
\label{alg:sweeping-setup}
\end{algorithm}

Once the preconditioner is set up, it can be applied using a straightforward
modification of Alg.~\ref{alg:naive-solve} which avoids 
redundant solves by combining the application of $L^{-1}$ and $D^{-1}$:

\begin{algorithm}
\DontPrintSemicolon
\tcp{Apply $L^{-1}$ and $D^{-1}$}
\For{$i=0,...,m-2$}{
  $u_i := T_i u_i$\;
  $u_{i+1} := u_{i+1} - J_{i+1,i} u_i$\;
}
$u_{m-1} := T_{m-1} u_{m-1}$\;
\tcp{Apply $L^{-T}$}
\For{$i=m-2,...,0$}{
  $u_i := u_i - T_i(J_{i+1,i}^T u_{i+1})$\;
}
\caption{Application of the sweeping preconditioner. 
$O(\gamma N \log N)$ work is required.}
\label{alg:sweeping-apply}
\end{algorithm}
Given that the preconditioner can be set up with $O(\gamma^2 N^{4/3})$ work, 
and applied with $O(\gamma N \log N)$ work, it provides a near-linear scheme
for solving Helmholtz equations if only $O(1)$ iterations are required for 
convergence. The experiments in this paper seem to indicate that, as long as 
the velocity model does not include a large cavity, the sweeping preconditioner 
indeed only requires $O(1)$ iterations.


Though this paper is focused on the parallel solution of Helmholtz equations,
which are the time-harmonic form of acoustic wave equations, Tsuji et al.\ have
shown that the moving PML sweeping preconditioner is equally effective
for time-harmonic Maxwell's equations~\cite{Tsuji-JCP,Tsuji-JFM}, and we believe
that the same will hold true for time-harmonic linear elasticity. The rest
of the paper will be presented in the context of the Helmholtz equation, but
we emphasize that the parallelization should carry over to more general wave 
equations in a conceptually trivial way. 

\subsection{Related work}
A domain decomposition variant of the sweeping preconditioner was recently
introduced~\cite{Stolk-sweeping} which results in fast convergence rates, 
albeit at the expense of requiring PML padding on both sides of each subdomain.
Recalling our previous analysis with respect to the PML size, 
$\gamma$, the memory usage of the preconditioner scales linearly with the 
PML size, while the setup cost scales quadratically. Thus, the domain 
decomposition approach should be expected to use twice as much memory and 
require four times as much work for the setup phase. 
On the other hand, doubling the subdomain sizes allows for more parallelism 
in both the setup and solve phases, and less sweeps seem to be required.

Another closely related method is the 
{\em Analytic ILU factorization}~\cite{GanderNataf-AILU}. 
Like the sweeping preconditioner, it uses local approximations of the Schur 
complements of the block $LDL^T$ factorization of the Helmholtz matrix 
represented in block tridiagonal form.
There are two crucial differences between the two methods:
\begin{itemize}
\item Roughly speaking, AILU can be viewed as using Absorbing Boundary 
      Conditions (ABC's)~\cite{EngquistMajda-ABC} instead of PML when forming 
      approximate subdomain auxiliary problems. 
      While ABC's result in strictly 2D local subproblems,
      versus the {\em quasi-2D} subdomain problems which result from using PML, 
      they are well-known to be less effective approximations of the Sommerfeld
      radiation condition (and thus the local Schur complement approximations 
      are less effective). The sweeping preconditioner handles its non-trivial
      subdomain factorizations via a multifrontal algorithm.
\item Rather than preconditioning with an approximate $LDL^T$ factorization
      of the original Helmholtz operator, the sweeping preconditioner sets up
      an approximate factorization of a {\em slightly damped} Helmholtz operator
      in order to mitigate the dispersion error which would result from the 
      Schur complement approximations.
\end{itemize}
These two improvements are responsible for transitioning from the $O(\omega)$ 
iterations required with AILU to the $O(1)$ iterations needed with the 
sweeping preconditioner (for problems without large cavities). 

Two other iterative methods warrant mentioning: the two-grid shifted-Laplacian 
approach of \cite{Calandra-twogrid} and the multilevel-ILU approach of 
\cite{BollhoeferGroteSchenk-ILU}. Though both require $O(\omega)$ iterations 
for convergence, they have very modest memory requirements.
In particular, \cite{Calandra-twogrid} demonstrates that, with a properly 
tuned two-grid approach, large-scale heterogeneous 3D problems can be solved 
with impressive timings.


There has also been a recent effort to extend the fast-direct methods presented 
in \cite{Xia-HSS} from definite elliptic problems into the realm of 
low-to-moderate frequency time-harmonic wave 
equations~\cite{Wang-geo,Wang-sisc}. 
In particular, their experiments (see Table 3 of \cite{Wang-geo}) suggest an 
asymptotic complexity of roughly $O(N^{1.8})$, which is a noticeable improvement
over the $O(N^2)$ complexity of traditional 3D sparse-direct methods.

\section{Parallel sweeping preconditioner}
\label{sec:parallel-sweeping}

The setup and application stages of the sweeping preconditioner (Algs.\ 
\ref{alg:sweeping-setup} and \ref{alg:sweeping-apply}) essentially consist of 
$m=O(n/\gamma)$ multifrontal factorizations and solves, respectively. 
The most important detail is that {\em the subdomain factorizations 
can be performed in parallel, while the subdomain solves must happen 
sequentially}.\footnote{While 
it is tempting to try to expose more parallelism with techniques like cyclic
reduction (which is a special case of a multifrontal algorithm), their
straightforward application destroys the Schur complement properties that we 
exploit for our fast algorithm.} When we also consider that each subdomain
factorization requires $O(\gamma^3 n^3)$ work, while subdomain solves only 
require $O(\gamma^2 n^2 \log n)$ work, we see that, relative to the subdomain 
factorizations, subdomain solves must extract a factor of $m=O(n/\gamma)$ more 
parallelism from a factor of $O(\gamma n/\log n)$ less operations. 
We thus have a strong hint that, unless the subdomain solves are carefully 
handled, they will be the limiting factor in the scalability of the sweeping 
preconditioner.

\subsection{Parallel multifrontal algorithms}
\label{subsection:parallel-multifrontal}
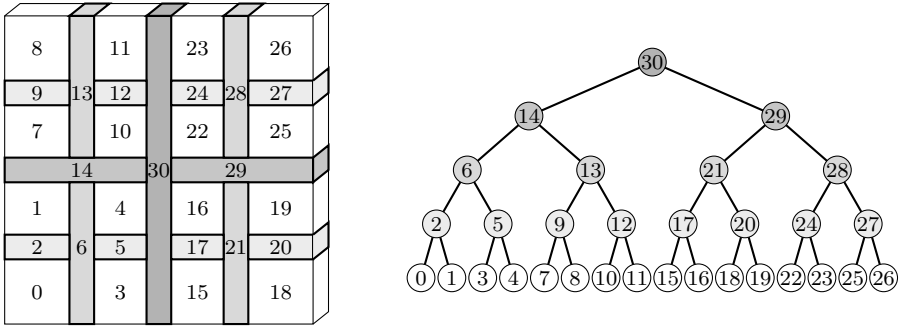
\begin{figure}
\centering
\begin{tikzpicture}[scale=4.1]
\tikzstyle{every node}=[font=\footnotesize]

\draw (0,0) rectangle (1,1);
\draw (0,1) -- (0.05,1.04) -- (1.05,1.04) -- (1,1);
\draw (1,0) -- (1.05,0.04) -- (1.05,1.04);

\draw[thick,fill=gray!60!] (0.46,0) rectangle (0.54,1); 
\draw (0.5,0.5) node{$30$};
\draw[thick,fill=gray!60!] 
  (0.46,1) -- (0.51,1.04) -- (0.59,1.04) -- (0.54,1) -- cycle;

\draw[thick,fill=gray!45!] (0,0.46) rectangle (0.46,0.54); 
\draw (0.25,0.5) node {$14$};
\draw[thick,fill=gray!45!] (0.54,0.46) rectangle (1,0.54); 
\draw (0.75,0.5) node {$29$};
\draw[thick,fill=gray!45!]
  (1,0.54) -- (1.05,0.58) -- (1.05,0.50) -- (1,0.46) -- cycle;

\draw[thick,fill=gray!30!] (0.21,0) rectangle (0.29,0.46); 
\draw (0.25,0.25) node {$6$};
\draw[thick,fill=gray!30!] (0.71,0) rectangle (0.79,0.46); 
\draw (0.75,0.25) node {$21$};
\draw[thick,fill=gray!30!] (0.21,0.54) rectangle (0.29,1); 
\draw (0.25,0.75) node {$13$};
\draw[thick,fill=gray!30!] 
  (0.21,1) -- (0.26,1.04) -- (0.34,1.04) -- (0.29,1) -- cycle;
\draw[thick,fill=gray!30!] (0.71,0.54) rectangle (0.79,1); 
\draw (0.75,0.75) node {$28$};
\draw[thick,fill=gray!30!]
  (0.71,1) -- (0.76,1.04) -- (0.84,1.04) -- (0.79,1) -- cycle;

\draw[thick,fill=gray!15!] (0,0.21) rectangle (0.21,0.29); 
\draw (0.105,0.25) node {$2$};
\draw[thick,fill=gray!15!] (0.29,0.21) rectangle (0.46,0.29); 
\draw (0.375,0.25) node {$5$};
\draw[thick,fill=gray!15!] (0,0.71) rectangle (0.21,0.79); 
\draw (0.105,0.75) node {$9$};
\draw[thick,fill=gray!15!] (0.29,0.71) rectangle (0.46,0.79); 
\draw (0.375,0.75) node {$12$};
\draw[thick,fill=gray!15!] (0.54,0.21) rectangle (0.71,0.29); 
\draw (0.625,0.25) node {$17$};
\draw[thick,fill=gray!15!] (0.79,0.21) rectangle (1,0.29); 
\draw (0.895,0.25) node {$20$};
\draw[thick,fill=gray!15!]
  (1,0.29) -- (1.05,0.33) -- (1.05,0.25) -- (1,0.21);
\draw[thick,fill=gray!15!] (0.54,0.71) rectangle (0.71,0.79); 
\draw (0.625,0.75) node {$24$};
\draw[thick,fill=gray!15!] (0.79,0.71) rectangle (1,0.79); 
\draw (0.895,0.75) node {$27$};
\draw[thick,fill=gray!15!]
  (1,0.79) -- (1.05,0.83) -- (1.05,0.75) -- (1,0.71);

\draw (0.105,0.105) node {$0$};
\draw (0.105,0.375) node {$1$};
\draw (0.375,0.105) node {$3$};
\draw (0.375,0.375) node {$4$};

\draw (0.105,0.625) node {$7$};
\draw (0.105,0.895) node {$8$};
\draw (0.375,0.625) node {$10$};
\draw (0.375,0.895) node {$11$};

\draw (0.625,0.105) node {$15$};
\draw (0.625,0.375) node {$16$};
\draw (0.895,0.105) node {$18$};
\draw (0.895,0.375) node {$19$};

\draw (0.625,0.625) node {$22$};
\draw (0.625,0.895) node {$23$};
\draw (0.895,0.625) node {$25$};
\draw (0.895,0.895) node {$26$};


\draw[thick] (1.35,0.15) -- (1.4,0.325) -- (1.45,0.15);
\draw[thick] (1.55,0.15) -- (1.6,0.325) -- (1.65,0.15);
\draw[thick] (1.75,0.15) -- (1.8,0.325) -- (1.85,0.15);
\draw[thick] (1.95,0.15) -- (2.0,0.325) -- (2.05,0.15);
\draw[thick] (2.15,0.15) -- (2.2,0.325) -- (2.25,0.15);
\draw[thick] (2.35,0.15) -- (2.4,0.325) -- (2.45,0.15);
\draw[thick] (2.55,0.15) -- (2.6,0.325) -- (2.65,0.15);
\draw[thick] (2.75,0.15) -- (2.8,0.325) -- (2.85,0.15);

\draw[thick] (1.4,0.325) -- (1.5,0.500) -- (1.6,0.325);
\draw[thick] (1.8,0.325) -- (1.9,0.500) -- (2.0,0.325);
\draw[thick] (2.2,0.325) -- (2.3,0.500) -- (2.4,0.325);
\draw[thick] (2.6,0.325) -- (2.7,0.500) -- (2.8,0.325);

\draw[thick] (1.5,0.500) -- (1.7,0.675) -- (1.9,0.500);
\draw[thick] (2.3,0.500) -- (2.5,0.675) -- (2.7,0.500);

\draw[thick] (1.7,0.675) -- (2.1,0.85) -- (2.5,0.675);

\draw[fill=white] (1.35,0.15) circle (0.045); \draw (1.35,0.15) node{$0$};
\draw[fill=white] (1.45,0.15) circle (0.045); \draw (1.45,0.15) node{$1$};
\draw[fill=white] (1.55,0.15) circle (0.045); \draw (1.55,0.15) node{$3$};
\draw[fill=white] (1.65,0.15) circle (0.045); \draw (1.65,0.15) node{$4$};
\draw[fill=white] (1.75,0.15) circle (0.045); \draw (1.75,0.15) node{$7$};
\draw[fill=white] (1.85,0.15) circle (0.045); \draw (1.85,0.15) node{$8$};
\draw[fill=white] (1.95,0.15) circle (0.045); \draw (1.95,0.15) node{$10$};
\draw[fill=white] (2.05,0.15) circle (0.045); \draw (2.05,0.15) node{$11$};
\draw[fill=white] (2.15,0.15) circle (0.045); \draw (2.15,0.15) node{$15$};
\draw[fill=white] (2.25,0.15) circle (0.045); \draw (2.25,0.15) node{$16$};
\draw[fill=white] (2.35,0.15) circle (0.045); \draw (2.35,0.15) node{$18$};
\draw[fill=white] (2.45,0.15) circle (0.045); \draw (2.45,0.15) node{$19$};
\draw[fill=white] (2.55,0.15) circle (0.045); \draw (2.55,0.15) node{$22$};
\draw[fill=white] (2.65,0.15) circle (0.045); \draw (2.65,0.15) node{$23$};
\draw[fill=white] (2.75,0.15) circle (0.045); \draw (2.75,0.15) node{$25$};
\draw[fill=white] (2.85,0.15) circle (0.045); \draw (2.85,0.15) node{$26$};

\draw[fill=gray!15!] (1.4,0.325) circle (0.045); \draw (1.4,0.325) node{$2$};
\draw[fill=gray!15!] (1.6,0.325) circle (0.045); \draw (1.6,0.325) node{$5$};
\draw[fill=gray!15!] (1.8,0.325) circle (0.045); \draw (1.8,0.325) node{$9$};
\draw[fill=gray!15!] (2.0,0.325) circle (0.045); \draw (2.0,0.325) node{$12$};
\draw[fill=gray!15!] (2.2,0.325) circle (0.045); \draw (2.2,0.325) node{$17$};
\draw[fill=gray!15!] (2.4,0.325) circle (0.045); \draw (2.4,0.325) node{$20$};
\draw[fill=gray!15!] (2.6,0.325) circle (0.045); \draw (2.6,0.325) node{$24$};
\draw[fill=gray!15!] (2.8,0.325) circle (0.045); \draw (2.8,0.325) node{$27$};

\draw[fill=gray!30!] (1.5,0.500) circle (0.045); \draw (1.5,0.500) node {$6$};
\draw[fill=gray!30!] (1.9,0.500) circle (0.045); \draw (1.9,0.500) node {$13$};
\draw[fill=gray!30!] (2.3,0.500) circle (0.045); \draw (2.3,0.500) node {$21$};
\draw[fill=gray!30!] (2.7,0.500) circle (0.045); \draw (2.7,0.500) node {$28$};

\draw[fill=gray!45!] (1.7,0.675) circle (0.045); \draw (1.7,0.675) node {$14$};
\draw[fill=gray!45!] (2.5,0.675) circle (0.045); \draw (2.5,0.675) node {$29$};

\draw[fill=gray!60!] (2.1,0.85) circle (0.045); \draw (2.1,0.85) node {$30$};

\end{tikzpicture}
\caption{A separator-based supernodal elimination tree (right) over a 
  quasi-2D subdomain (left).}
\label{fig:sep-tree}
\end{figure}

\begin{figure}
\centering
\begin{tikzpicture}[scale=4.1]
\tikzstyle{every node}=[font=\footnotesize]

\draw[thick] (0.1,0.15) -- (0.2,0.325) -- (0.3,0.15);
\draw[thick] (0.5,0.15) -- (0.6,0.325) -- (0.7,0.15);
\draw[thick] (0.9,0.15) -- (1.0,0.325) -- (1.1,0.15);
\draw[thick] (1.3,0.15) -- (1.4,0.325) -- (1.5,0.15);
\draw[thick] (1.7,0.15) -- (1.8,0.325) -- (1.9,0.15);
\draw[thick] (2.1,0.15) -- (2.2,0.325) -- (2.3,0.15);
\draw[thick] (2.5,0.15) -- (2.6,0.325) -- (2.7,0.15);
\draw[thick] (2.9,0.15) -- (3.0,0.325) -- (3.1,0.15);

\draw[thick] (0.2,0.325) -- (0.4,0.500) -- (0.6,0.325);
\draw[thick] (1.0,0.325) -- (1.2,0.500) -- (1.4,0.325);
\draw[thick] (1.8,0.325) -- (2.0,0.500) -- (2.2,0.325);
\draw[thick] (2.6,0.325) -- (2.8,0.500) -- (3.0,0.325);

\draw[thick] (0.4,0.500) -- (0.8,0.675) -- (1.2,0.500);
\draw[thick] (2.0,0.500) -- (2.4,0.675) -- (2.8,0.500);

\draw[thick] (0.8,0.675) -- (1.6,0.85) -- (2.4,0.675);

\draw[fill=white] (0.1,0.15) ellipse (0.09 and 0.045); 
\draw (0.1,0.15) node{$000$};
\draw[fill=white] (0.3,0.15) ellipse (0.09 and 0.045); 
\draw (0.3,0.15) node{$000$};
\draw[fill=white] (0.5,0.15) ellipse (0.09 and 0.045); 
\draw (0.5,0.15) node{$001$};
\draw[fill=white] (0.7,0.15) ellipse (0.09 and 0.045); 
\draw (0.7,0.15) node{$001$};
\draw[fill=white] (0.9,0.15) ellipse (0.09 and 0.045); 
\draw (0.9,0.15) node{$010$};
\draw[fill=white] (1.1,0.15) ellipse (0.09 and 0.045); 
\draw (1.1,0.15) node{$010$};
\draw[fill=white] (1.3,0.15) ellipse (0.09 and 0.045); 
\draw (1.3,0.15) node{$011$};
\draw[fill=white] (1.5,0.15) ellipse (0.09 and 0.045); 
\draw (1.5,0.15) node{$011$};
\draw[fill=white] (1.7,0.15) ellipse (0.09 and 0.045); 
\draw (1.7,0.15) node{$100$};
\draw[fill=white] (1.9,0.15) ellipse (0.09 and 0.045); 
\draw (1.9,0.15) node{$100$};
\draw[fill=white] (2.1,0.15) ellipse (0.09 and 0.045); 
\draw (2.1,0.15) node{$101$};
\draw[fill=white] (2.3,0.15) ellipse (0.09 and 0.045); 
\draw (2.3,0.15) node{$101$};
\draw[fill=white] (2.5,0.15) ellipse (0.09 and 0.045); 
\draw (2.5,0.15) node{$110$};
\draw[fill=white] (2.7,0.15) ellipse (0.09 and 0.045); 
\draw (2.7,0.15) node{$110$};
\draw[fill=white] (2.9,0.15) ellipse (0.09 and 0.045); 
\draw (2.9,0.15) node{$111$};
\draw[fill=white] (3.1,0.15) ellipse (0.09 and 0.045); 
\draw (3.1,0.15) node{$111$};

\draw[fill=gray!15!] (0.2,0.325) ellipse (0.09 and 0.045); 
\draw (0.2,0.325) node{$000$};
\draw[fill=gray!15!] (0.6,0.325) ellipse (0.09 and 0.045); 
\draw (0.6,0.325) node{$001$};
\draw[fill=gray!15!] (1.0,0.325) ellipse (0.09 and 0.045); 
\draw (1.0,0.325) node{$010$};
\draw[fill=gray!15!] (1.4,0.325) ellipse (0.09 and 0.045); 
\draw (1.4,0.325) node{$011$};
\draw[fill=gray!15!] (1.8,0.325) ellipse (0.09 and 0.045); 
\draw (1.8,0.325) node{$100$};
\draw[fill=gray!15!] (2.2,0.325) ellipse (0.09 and 0.045); 
\draw (2.2,0.325) node{$101$};
\draw[fill=gray!15!] (2.6,0.325) ellipse (0.09 and 0.045); 
\draw (2.6,0.325) node{$110$};
\draw[fill=gray!15!] (3.0,0.325) ellipse (0.09 and 0.045); 
\draw (3.0,0.325) node{$111$};

\draw[fill=gray!30!] (0.4,0.500) ellipse (0.09 and 0.045); 
\draw (0.4,0.500) node {$00*$};
\draw[fill=gray!30!] (1.2,0.500) ellipse (0.09 and 0.045); 
\draw (1.2,0.500) node {$01*$};
\draw[fill=gray!30!] (2.0,0.500) ellipse (0.09 and 0.045); 
\draw (2.0,0.500) node {$10*$};
\draw[fill=gray!30!] (2.8,0.500) ellipse (0.09 and 0.045); 
\draw (2.8,0.500) node {$11*$};

\draw[fill=gray!45!] (0.8,0.675) ellipse (0.09 and 0.045); 
\draw (0.8,0.675) node {$0**$};
\draw[fill=gray!45!] (2.4,0.675) ellipse (0.09 and 0.045); 
\draw (2.4,0.675) node {$1**$};

\draw[fill=gray!60!] (1.6,0.85) ellipse (0.09 and 0.045); 
\draw (1.6,0.85) node {$***$};

\end{tikzpicture}
\caption{Overlay of the process ranks (in binary) of the owning subteams of each
supernode from the elimination tree in Fig.\ \ref{fig:sep-tree} when the tree 
is assigned to eight processes using a subtree-to-subteam mapping; a `*' is 
used to denote both 0 and 1, so that `$00*$' represents processes 0 and 1, 
`$01*$' represents processes 2 and 3, and `$***$' represents all eight 
processes.}
\label{fig:subteam}
\end{figure}
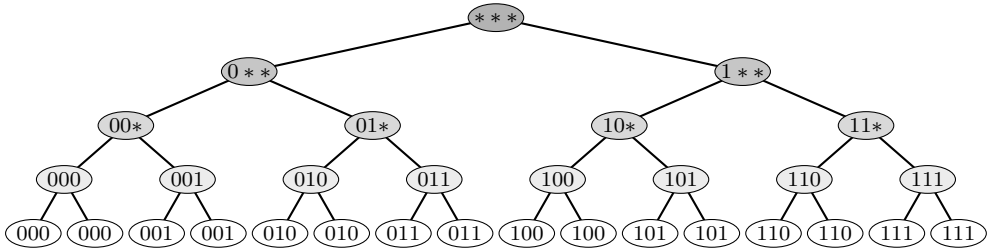

While a large number of techniques exist for parallelizing multifrontal
factorizations and triangular solves, we focus on parallelizations 
which combine subtree-to-subteam~\cite{GeorgeLiuNg-subtree} mappings of 
processes to the elimination tree~\cite{Schreiber-etree} that also make use of 
two-dimensional distributions of the frontal 
matrices~\cite{Schreiber-scalability}.\footnote{Cf.~\cite{Amestoy-MUMPS}, which
advocates for only distributing the root frontal matrix two-dimensionally and 
using a one-dimensional distribution for all other fronts.}
More specifically, we make use of supernodal~\cite{Ashcraft-progress} 
elimination trees defined through nested dissection
(see Figs.~\ref{fig:sep-tree} and \ref{fig:subteam}), which have been shown to 
result in highly scalable 
factorizations~\cite{GuptaKarypisKumar-scalable,GuptaKoricGeorge-massive} and 
moderately scalable triangular solutions~\cite{JoshiGuptaKarypisKumar-2d}. 


Roughly speaking, the analysis in \cite{JoshiGuptaKarypisKumar-2d} shows that, 
if $p_F$ processes are used in the multifrontal factorization of our quasi-2D 
subdomain problems, then we must have $\gamma n=\Omega(p_F^{1/2})$ in order to 
maintain constant efficiency as $p_F$ is increased; similarly, if $p_S$ 
processes are used in the multifrontal triangular solves for a subdomain, 
then we must have $\gamma n\approx \Omega(p_S)$
(where we use $\approx$ to denote that the equality holds within logarithmic 
factors). 
Since we can simultaneously factor the $m=O(n/\gamma)$ subdomain 
matrices, we denote the total number of processes as $p$ and set 
$p_S=p$ and $p_F=O(p/m)$; then the subdomain factorizations only require
that $n^3=\Omega(p/\gamma)$, while the subdomain solves have the 
much stronger constraint that $n \approx \Omega(p/\gamma)$. 
This last constraint should be considered unacceptable, as we have the 
conflicting requirement that $n^3 \approx O(p/\gamma)$ in 
order to store the factorizations in memory. It is therefore advantageous to 
consider more scalable alternatives to standard multifrontal triangular solves, 
even if they require additional computation.

\subsection{Selective inversion}
\label{subsection:selective-inversion}
The lackluster scalability of dense triangular solves is well known and
a scheme known as {\em selective inversion} was introduced 
in~\cite{Raghavan-invert} and implemented in~\cite{Raghavan-dscpack} 
specifically to avoid the issue; the approach is 
characterized by directly inverting every distributed dense triangular matrix 
which would have been solved against in a normal multifrontal triangular solve. Thus, after performing selective inversion, every parallel dense triangular 
solve can be translated into a parallel dense triangular matrix-vector multiply.

Suppose that we have paused a multifrontal $LDL^T$ factorization just before
processing a particular front, $F$, which corresponds to some supernode, 
$\mathcal{S}$. Then all of the fronts for the descendants of $\mathcal{S}$ 
have already been handled, and $F$ can be partitioned as
\begin{equation}
  F = \left(\begin{array}{cc} F_{TL} & \star \\ 
                              F_{BL} & F_{BR}\end{array}\right),
\end{equation}
where $F_{TL}$ holds the Schur complement of supernode $\mathcal{S}$ with 
respect to all of its descendants, $F_{BL}$ 
represents the coupling of $\mathcal{S}$ and its descendants to $\mathcal{S}$'s
ancestors, and $F_{BR}$ holds the Schur complement updates from 
the descendants of $\mathcal{S}$ for the ancestors of $\mathcal{S}$. Using 
hats to denote input states, e.g., $\hat F_{TL}$ to denote the input state
of $F_{TL}$, the first step in processing the frontal matrix $F$ is to 
overwrite $F_{TL}$ with its $LDL^T$ factorization, which is to say that 
$\hat F_{TL}$ is overwritten with the strictly lower portion of a unit lower 
triangular matrix $L_F$ and a diagonal matrix $D_F$ such that 
$\hat F_{TL} = L_F D_F L_F^T$. 

The partial factorization of $F$ can then be completed via the following steps:
\begin{enumerate}
\item Solve against $L_F^T$ to form $F_{BL} := F_{BL} L_F^{-T}$.
\item Form the temporary copy $Z_{BL} := F_{BL}$.
\item Finalize the coupling matrix as $F_{BL} := F_{BL} D_F^{-1}$.
\item Finalize the update matrix as 
$F_{BR} := F_{BR} - \hat F_{BL} \hat F_{TL}^{-1} \hat F_{BL}^T 
         = F_{BR} - Z_{BL} F_{BL}^T$.
\end{enumerate}
After adding $F_{BR}$ onto the parent frontal matrix, only $F_{TL}$ and $F_{BL}$
are needed in order to perform a multifrontal solve. For instance, applying
$L^{-1}$ to some vector $x$ proceeds up the elimination tree 
(starting from the leaves) in a manner similar to the factorization;
after handling all of the work for the descendants of some supernode 
$\mathcal{S}$, only a few dense linear algebra operations with $\mathcal{S}$'s 
corresponding frontal matrix, say $F$, are required.
Denoting the portion of $x$ corresponding to the degrees of freedom in 
supernode $\mathcal{S}$ by $x_{\mathcal{S}}$, we must perform:
\begin{enumerate}
\item $x_{\mathcal{S}} := L_F^{-1} x_{\mathcal{S}}$
\item $x_U \equiv -F_{BL} x_{\mathcal{S}}$
\item Add $x_U$ onto the entries of $x$ corresponding to the parent supernode.
\end{enumerate}
The key insight of selective inversion is that, 
{\em if we invert each distributed dense unit lower triangular matrix $L_F$ in 
place, all of the parallel dense triangular solves in a multifrontal 
triangular solve are replaced by parallel dense matrix-vector multiplies.} It 
is also observed in \cite{Raghavan-invert} that the work required for the 
selective inversion is typically only a modest percentage of the work required 
for the multifrontal factorization, and that the overhead of the selective 
inversion is easily recouped if there are several right-hand sides to solve 
against.

Since each application of the sweeping preconditioner requires two multifrontal
solves for each of the $m=O(n/\gamma)$ subdomains, which are relatively small 
and likely distributed over a large number of processes, selective inversion 
will be shown to yield a very large performance improvement.
We also note that, while it is widely believed that direct inversion is 
numerically unstable, in \cite{DruinskyToledo-inverse} Druinsky and Toledo 
provide a review of (apparently obscure) results dating back to Wilkinson 
(in \cite{Wilkinson-rounding}) which show that $x := \mathrm{inv}(A)\!*\!b$ is
as accurate as a backwards stable solve if reasonable assumptions are met
on the accuracy of $\mathrm{inv}(A)$. Since $\mathrm{inv}(A)\!*\!b$ is argued to
be more accurate when the columns of $\mathrm{inv}(A)$ have been computed 
with a backwards-stable solver, and both $\mathrm{inv}(F_{TL})$ and 
$\mathrm{inv}(F_{TL}^T)$ must be applied after selective inversion, it might
be worthwhile to modify selective inversion to compute and store two different 
inverses of each $F_{TL}$: one by columns and one by rows. 

\subsection{Global vector distributions}
The goal of this subsection is to determine an appropriate scheme for
distributing global vectors, i.e., ones representing a function over the 
entire domain (as opposed to only over a panel). And while the
factorizations themselves may have occurred on subteams of $O(p/m)$ processes 
each, in order to make use of all available processes for every subdomain solve,
at this point we assume that each auxiliary problem's frontal tree has 
been selectively inverted and is distributed using a subtree-to-subteam mapping 
(recall Fig.~\ref{fig:subteam}) over the entire set of $p$ 
processes.\footnote{In cases where the available solve parallelism has been 
exhausted but the problem cannot be solved on less processes due to memory 
constraints, it would be preferable to solve with subdomains stored on subsets 
of processes.}

Since a subtree-to-subteam mapping will assign each supernode of an 
auxiliary problem to a team of processes, and each panel of the original
3D domain is by construction a subset of the domain of an auxiliary 
problem, it is straightforward to extend the supernodal subteam assignments to 
the full domain. We should then be able to distribute global vectors so
that no communication is required for readying panel subvectors for 
subdomain solves (via extension by zero for interior panels, and trivially for 
the first panel). Since our nested dissection process does not partition in 
the shallow dimension of quasi-2D subdomains (see Fig.~\ref{fig:sep-tree}), 
extending a vector defined over a panel of the original domain onto the 
PML-padded auxiliary domain simply requires individually extending each 
supernodal subvector by zero in the $x_3$ direction.

Consider an element-wise two-dimensional cyclic distribution 
\cite{Poulson-elemental} of a frontal matrix $F$ over $q$ processes using an 
$r \times c$ process grid, where $r$ and $c$ are $O(\sqrt{q})$. Then the 
$(i,j)$ entry will be stored by the process in the $(i \bmod r,j \bmod c)$ 
position in the process grid. 
Using the notation from \cite{Poulson-elemental}, this distributed front would 
be denoted as $F[M_C,M_R]$, while its top-left quadrant would be referred to 
as $F_{TL}[M_C,M_R]$ (see Fig.~\ref{fig:mc_mr} for a depiction of an $[M_C,M_R]$ matrix distribution). 

\begin{figure}
\centering
$
\left(\begin{array}{cccccccc}
 0 & 2 & 4 & 0 & 2 & 4 & 0 \\
 1 & 3 & 5 & 1 & 3 & 5 & 1 \\
 0 & 2 & 4 & 0 & 2 & 4 & 0 \\
 1 & 3 & 5 & 1 & 3 & 5 & 1 \\
 0 & 2 & 4 & 0 & 2 & 4 & 0 \\
 1 & 3 & 5 & 1 & 3 & 5 & 1 \\
 0 & 2 & 4 & 0 & 2 & 4 & 0 
\end{array}\right)\;\;\;\;\;\;\;
\begin{array}{ccccc}
0 & - & 2 & - & 4 \\
| &   & | &   & | \\
1 & - & 3 & - & 5
\end{array}
$
\caption{Overlay of the owning process ranks of an 7 $\times$ 7 matrix 
distributed over a 2 $\times$ 3 process grid in the $[M_C,M_R]$ distribution, 
where $M_C$ assigns row $i$ to process row $i \bmod 2$, and $M_R$ assigns column
$j$ to process column $i \bmod 3$ (left). The process grid is shown on the 
right.}
\label{fig:mc_mr}
\end{figure}

Loosely speaking, $F[X,Y]$ means that each column
of $F$ is distributed using the scheme denoted by $X$, and each row
is distributed using the scheme denoted by $Y$. For the element-wise 
two-dimensional distribution used for $F$, $[M_C,M_R]$, we have that the 
columns of $F$ are distributed like Matrix Columns ($M_C$), and the rows 
of $F$ are distributed like Matrix Rows ($M_R$). While this notation may seem
vapid when only working with a single distributed matrix, it is useful when
working with products of distributed matrices. For instance, if a `$\star$' is 
used to represent rows/columns being redundantly stored (i.e., not distributed),
then the result of every process multiplying its local submatrix of 
$A[X,\star]$ with its local submatrix of $B[\star,Y]$ forms a distributed
matrix $C[X,Y] = (AB)[X,Y] = A[X,\star]\, B[\star,Y]$, where the last 
expression refers to the local multiplication process.

We can now decide on a distribution for each supernodal subvector, say 
$x_{\mathcal{S}}$, based on the criteria that it should be fast to 
form $F_{TL} x_{\mathcal{S}}$ and $F_{TL}^T x_{\mathcal{S}}$ in the same 
distribution as $x_{\mathcal{S}}$, given that $F_{TL}$ is distributed as 
$F_{TL}[M_C,M_R]$. Suppose that we define a Column-major Vector distribution 
($V_C$) of $x_{\mathcal{S}}$, say $x_{\mathcal{S}}[V_C,\star]$, to 
mean that entry $i$ is owned by process $i \bmod q$, which corresponds 
to position $(i \bmod r,\lfloor i/r \rfloor \bmod c)$ in the process grid 
(if the grid is constructed with a column-major ordering of the process ranks; 
 see the left side of Fig.~\ref{fig:vector}).
Then a call to \verb!MPI_Allgather!~\cite{Dongarra-mpi} within each row of the 
process grid would allow for each process to collect all of the data necessary 
to form $x_{\mathcal{S}}[M_C,\star]$, as for any process row index 
$s \in \{0,1,...,r-1\}$,
\begin{equation}\label{vc-to-mc}
  \{ i \in \mathbb{N}_0 : i \bmod r = s \} = 
  \bigcup_{t=0}^{c-1} \{ i \in \mathbb{N}_0 : i \bmod q = s+tr \}.
\end{equation}
See the left side of Fig.~\ref{fig:partial_matrix} for an example of an 
$[M_C,\star]$ distribution of a $7 \times 3$ matrix.

\begin{figure}
\centering
$
\left(\begin{array}{c}
 0 \\
 1 \\
 2 \\
 3 \\
 4 \\
 5 \\
 0  
\end{array}\right),\;\;\;\;\;\;
\left(\begin{array}{c}
 0 \\
 2 \\
 4 \\
 1 \\
 3 \\
 5 \\
 0 
\end{array}\right)
$
\caption{Overlay of the owning process ranks of a vector of height 7 
distributed over a 2 $\times$ 3 process grid in the $[V_C,\star]$ vector 
distribution (left) and the $[V_R,\star]$ vector distribution (right).}
\label{fig:vector}
\end{figure}

\begin{figure}
\centering
$
\left(\begin{array}{c}
 \{0,2,4\} \\
 \{1,3,5\} \\
 \{0,2,4\} \\
 \{1,3,5\} \\
 \{0,2,4\} \\
 \{1,3,5\} \\
 \{0,2,4\} 
\end{array}\right),\;\;\;\;\;\;
\left(\begin{array}{c}
 \{0,1\} \\
 \{2,3\} \\
 \{4,5\} \\
 \{0,1\} \\
 \{2,3\} \\
 \{4,5\} \\
 \{0,1\} 
\end{array}\right)
$
\caption{Overlay of the owning process ranks of a vector of height 7 
distributed over a 2 $\times$ 3 process grid in the $[M_C,\star]$ 
distribution (left) and the $[M_R,\star]$ distribution (right).}
\label{fig:partial_matrix}
\end{figure}

Similarly, if $x_{\mathcal{S}}$ was distributed with a Row-major Vector 
distribution ($V_R$), as shown on the right side of Fig.~\ref{fig:vector}, 
say $x_{\mathcal{S}}[V_R,\star]$, so that 
entry $i$ is owned by the process in position 
$(\lfloor i/c \rfloor \bmod r,i \bmod c)$ of the process grid, then a call to
\verb!MPI_Allgather! within each column of the process grid would provide 
each process with the data necessary to form $x_{\mathcal{S}}[M_R,\star]$.
Under reasonable assumptions, both of these redistributions can be shown to 
have per-process communication volume lower bounds of $\Omega(n/\sqrt{p})$ 
(if $F_{TL}$ is $n \times n$) and latency lower bounds of 
$\Omega(\log_2(\sqrt{p}))$~\cite{Chan-collective}. 
We also note that translating between $x_{\mathcal{S}}[V_C,\star]$ and 
$x_{\mathcal{S}}[V_R,\star]$ simply requires permuting which process owns each 
local subvector, so the communication volume would be $O(n/p)$, 
while the latency cost is $O(1)$.

We have thus described efficient techniques for redistributing 
$x_{\mathcal{S}}[V_C,\star]$ to both the $x_{\mathcal{S}}[M_R,\star]$ and 
$x_{\mathcal{S}}[M_C,\star]$ distributions, which are the first steps 
for our parallel algorithms for forming $F_{TL} x_{\mathcal{S}}$ and 
$F_{TL}^T x_{\mathcal{S}}$, respectively: Denoting the distributed result of 
each process in process column $t \in \{0,1,...,c-1\}$ multiplying its local 
submatrix of $F_{TL}[M_C,M_R]$ by its local subvector of 
$x_{\mathcal{S}}[M_R,\star]$ as $z^{(t)}[M_C,\star]$, it holds that 
$(F_{TL} x_{\mathcal{S}})[M_C,\star] = \sum_{t=0}^{c-1} z^{(t)}[M_C,\star]$.
Since Eq.\ \eqref{vc-to-mc} also implies that each process's local data 
from a $[V_C,\star]$ distribution is a subset of its local data from a 
$[M_C,\star]$ distribution, a simultaneous summation and scattering of 
$\{z^{(t)}[M_C,\star]\}_{t=0}^{c-1}$ within process rows, perhaps via 
\verb!MPI_Reduce_scatter! or \verb!MPI_Reduce_scatter_block!, yields the 
desired result, $(F_{TL} x_{\mathcal{S}})[V_C,\star]$. An analogous process
with $(F_{TL}[M_C,M_R])^T=F_{TL}^T[M_R,M_C]$ and $x_{\mathcal{S}}[M_C,\star]$ 
yields $(F_{TL}^T x_{\mathcal{S}})[V_R,\star]$, which can then be cheaply 
permuted to form $(F_{TL}^T x_{\mathcal{S}})[V_C,\star]$. Both calls to 
\verb!MPI_Reduce_scatter_block! can be shown to have the same communication 
lower bounds as the previously discussed \verb!MPI_Allgather! 
calls~\cite{Chan-collective}.

As discussed at the beginning of this section, defining the distribution of 
each supernodal subvector specifies a distribution for a global vector,
say $[\mathcal{G},\star]$. While the $[V_C,\star]$ distribution shown in  
the left half of Fig.~\ref{fig:vector} simply assigns entry $i$ of a 
supernodal subvector $x_{\mathcal{S}}$, distributed over $q$ processes, to 
process $i \bmod q$, we can instead choose an alignment 
parameter, $\sigma$, where $0 \le \sigma < q$, and assign entry $i$ to 
process $(i + \sigma) \bmod q$. If we simply set $\sigma=0$ for every 
supernode, as the discussion at the beginning of this subsection implied, then 
at most $O(\gamma n)$ processes will store data for the root separator 
supernodes of a global vector, as each root separator only has $O(\gamma n)$ 
degrees of freedom by construction.
However, there are $m=O(n/\gamma)$ root separators, so we can easily allow 
for up to $O(n^2)$ processes to share the storage of a global vector if the 
alignments are carefully chosen. 
It is important to notice that the top-left quadrants of the 
frontal matrices for the root separators can each be distributed over 
$O(\gamma^2 n^2)$ processes, so $O(\gamma^2 n^2)$ processes can actively 
participate in the corresponding triangular matrix-vector multiplications.

\subsection{Parallel preconditioned GMRES(k)}
Since, by hypothesis, only $O(1)$ iterations of \gmresk~will be required for 
convergence with the sweeping preconditioner, a cursory inspection of 
Algorithm~\ref{alg:sweeping-apply} reveal that most of the work in a 
preconditioned iterative method, such as \gmresk, will be performed in the 
multifrontal solves during the preconditioner application, but a modest 
portion will also be spent in sparse matrix-vector multiplication with the 
discrete Helmholtz operator, $A$, and the off-diagonal blocks of the discrete 
artificially damped Helmholtz operator, 
$J$. It is thus important to parallelize the sparse matrix-vector multiplies, 
but it is not crucial that the scheme be optimal, and so we simply distribute
$A$ and $J$ in the same manner as vectors, i.e., with the $[\mathcal{G},\star]$ 
distribution derived from the auxiliary problems' frontal distributions. 


For performance reasons, it is beneficial to solve as many right-hand sides 
simultaneously as possible: both the communication latency and the costs of 
loading the local data from frontal and sparse matrices from main memory 
can be amortized over all of the right-hand sides. Another idea is to extend 
the so-called \verb!trsm! algorithm for triangular solves 
with many right-hand sides (i.e., more right-hand sides than processes), 
which is well-known in the field of dense linear 
algebra~\cite{Poulson-elemental}, into the realm of sparse-direct solvers via
the dense frontal triangular solves.
This approach was not pursued in this paper due to the modest storage space
available on Lonestar and is left for future work.
Another performance improvement might come from exploiting block variants of 
GMRES~\cite{Simoncini-block-gmres}, which can potentially lower the number of 
required iterations. 

\subsection{Clique}
In order to implement the previously discussed techniques for scalable 
multifrontal factorizations and solves (via selective inversion), an 
open-source distributed multifrontal solver named Clique was built on top of 
Elemental~\cite{Poulson-elemental}, a library for distributed-memory dense 
linear algebra.  In addition to being designed to support the techniques 
we discussed above: selective inversion, subtree-to-submesh mappings, and 
two-dimensional frontal matrix distributions, it was also written with 
a strong emphasis on {\em memory scalability}. This is because the 
sweeping preconditioner requires large numbers of factorizations 
of relatively small sparse matrices, and so it is crucial that the per-process
memory usage for each subdomain factorization decreases inversely with the 
total number of processes.

We note that Clique was designed specifically to provide a memory-scalable 
multifrontal implementation for our parallel sweeping preconditioner, and so
there is not yet support for pivoting. We plan to add a pivoted LU factorization
in the near future.

\subsection{Parallel Sweeping Preconditioner (PSP)}
Given the discussion in Section \ref{sec:parallel-sweeping}, it is most 
convenient to describe our prototype implementation of a parallel sweeping 
preconditioner based upon its deviations from our proposed approach. The 
primary difference is that there is not yet support for simultaneously 
factoring the subdomain auxiliary problems and then redistributing each frontal
tree to the entire set of processes. This will certainly lead to large
improvements in the scalability of the setup phase, but it is left for future 
work.

\section{Experimental results}
\label{sec:experimental}
Our experiments were performed on the Texas Advanced Computing Center
(TACC) machine, Lonestar, which is comprised of 1,888 compute nodes, 
each equipped with two hex-core 3.33 GHz processors and 24 GB of memory, 
which are connected with QDR InfiniBand using a fat-tree topology. 
Our tests launched eight MPI processes per node in order to provide each MPI 
process with 3 GB of memory.

\begin{figure}
\centering
\includegraphics[width=3.8in]{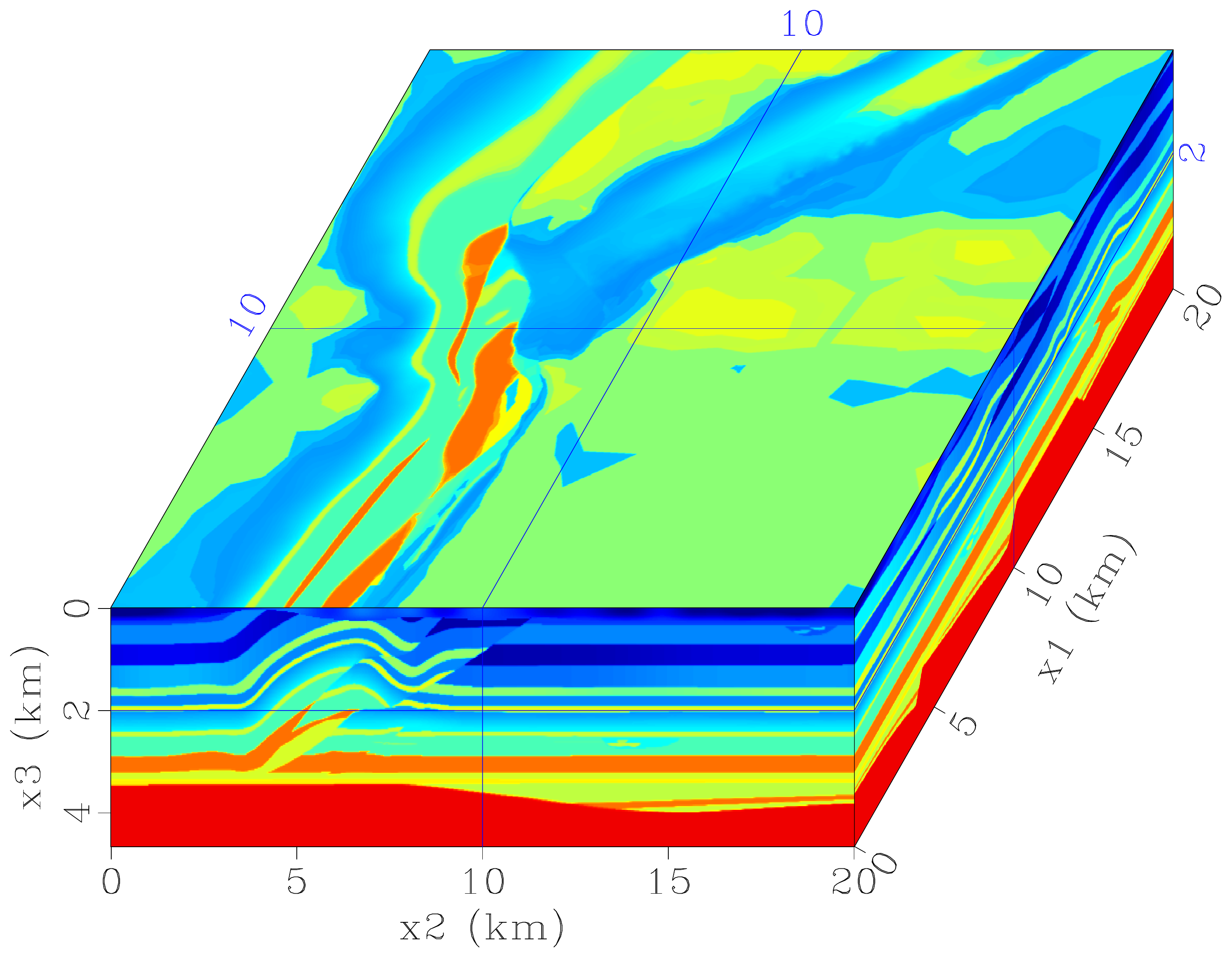}
\caption{Three cross-sections of the SEG/EAGE Overthrust velocity model, which 
represents an artificial 
$20\, \mathrm{km} \times 20\, \mathrm{km} \times 4.65\, \mathrm{km}$ domain, 
containing an {\em overthrust} fault, using samples every $25\, \mathrm{m}$. 
The result is an $801 \times 801 \times 187$ grid of wave speeds varying 
discontinuously between $2.179\, \mathrm{km/sec}$ (blue) and 
$6.000\, \mathrm{km/sec}$ (red).}
\label{fig:overthrust}
\end{figure}

Our experiments took place over five different 3D velocity models:
\begin{itemize}
\item A uniform background with a high-contrast barrier. The domain is 
      the unit cube and the wave speed is 1 except in 
      $[0,1] \times [0.25,0.3] \times [0,0.75]$, where it is $10^{10}$.
\item A wedge problem over the unit cube, where the wave speed is set to $2$ 
      if $Z \le 0.4+0.1x_2$, $1.5$ if otherwise $Z \le 0.8 - 0.2x_2$, and $3$ in
      all other cases.
\item A two-layer model defined over the unit cube, 
      where $c=4$ if $x_2<0.5$, and $c=1$ otherwise.
\item A waveguide over the unit cube:
      $c(\mathbf{x})=1.25(1-0.4 e^{-32 (|x_1-0.5|^2+|x_2-0.5|^2)})$.
\item The SEG/EAGE Overthrust model~\cite{Aminzadeh-overthrust}, 
      see Fig.~\ref{fig:overthrust}.
\end{itemize}

In all of the following experiments, the shortest wavelength of each model is
resolved with roughly ten grid points and the performance of the preconditioner
is tested using the following four forcing functions:
\begin{itemize}
\item a single {\em shot} centered at $\mathbf{x}_0$,
$f_0(\mathbf{x}) = n e^{-10 n \|\mathbf{x}-\mathbf{x}_0\|^2}$,
\item three shots,
$f_1(\mathbf{x}) = \sum_{i=0}^2 n e^{-10 n \|\mathbf{x}-\mathbf{x}_i\|^2}$,
\item a Gaussian beam centered at $\mathbf{x}_2$ in direction $\mathbf{d}$,
$f_2(\mathbf{x})=e^{i\omega \mathbf{x} \cdot \mathbf{d}} 
e^{-4\omega \|\mathbf{x}-\mathbf{x}_2\|^2}$, and 
\item a plane wave in direction $\mathbf{d}$, 
$f_3(\mathbf{x})=e^{i\omega \mathbf{x} \cdot \mathbf{d}}$,
\end{itemize}
where $\mathbf{x}_0=(0.5,0.5,0.1)$, $\mathbf{x}_1=(0.25,0.25,0.1)$, 
$\mathbf{x}_2=(0.75,0.75,0.5)$, and $\mathbf{d}=(1,1,-1)/\sqrt{3}$.
Note that, in the case of the Overthrust model, these source locations should 
be interpreted proportionally (e.g., an $x_3$ value of $0.1$ means a depth
which is $10\%$ of the total depth of the model). Due to the oscillatory nature
of the solution, we simply choose the zero vector as our initial guess in 
all experiments.

The first experiment was meant to test the convergence rate of the sweeping 
preconditioner over domains spanning 50 wavelengths in each direction
(resolved to ten points per wavelength) and each test made use of 256 nodes of 
Lonestar. During the course of the tests, it was noticed that a significant 
amount of care must be taken when setting the amplitude of the derivative of 
the PML takeoff function (i.e., the ``C'' variable in Eq.~(2.1) 
from~\cite{EngquistYing-PML}). For the sake of brevity, hereafter we refer to
this variable as the {\em PML amplitude}. A modest search was performed in
order to find a near-optimal value for the PML amplitude for each problem.
These values are reported in Table \ref{tbl:convergence} along with the number 
of iterations required for the relative residuals for all four forcing functions
to reduce to less than $10^{-5}$.

\begin{table}
\small
\centering
\begin{tabular}{r|c|c|c|c|}
\cline{2-5}
 & \multicolumn{4}{|c|}{velocity model} \\
\cline{2-5}
  & barrier & wedge & two-layers & waveguide  \\
\hline
\multicolumn{1}{|r|}{Hz}
  & 50 & 75 & 50 & 37.5 \\
\hline
\multicolumn{1}{|r|}{PML amplitude}
  & 3.0 & 4.0 & 4.0 & 2.0 \\
\hline
\multicolumn{1}{|r|}{iterations}
  & 28 & 49 & 48 & 52 \\
\hline
\end{tabular}
\caption{The number of iterations required for convergence for four model 
problems (with four forcing functions per model).
The grid sizes were $500^3$ and roughly 50 wavelengths were spanned in each 
direction. Five grid points were used for all PML discretizations, 
four planes were processed per panel, and the damping factors were all set to 
$7$.}
\label{tbl:convergence}
\end{table}

\begin{figure}
$
\begin{array}{cc}
\includegraphics[width=2.5in]{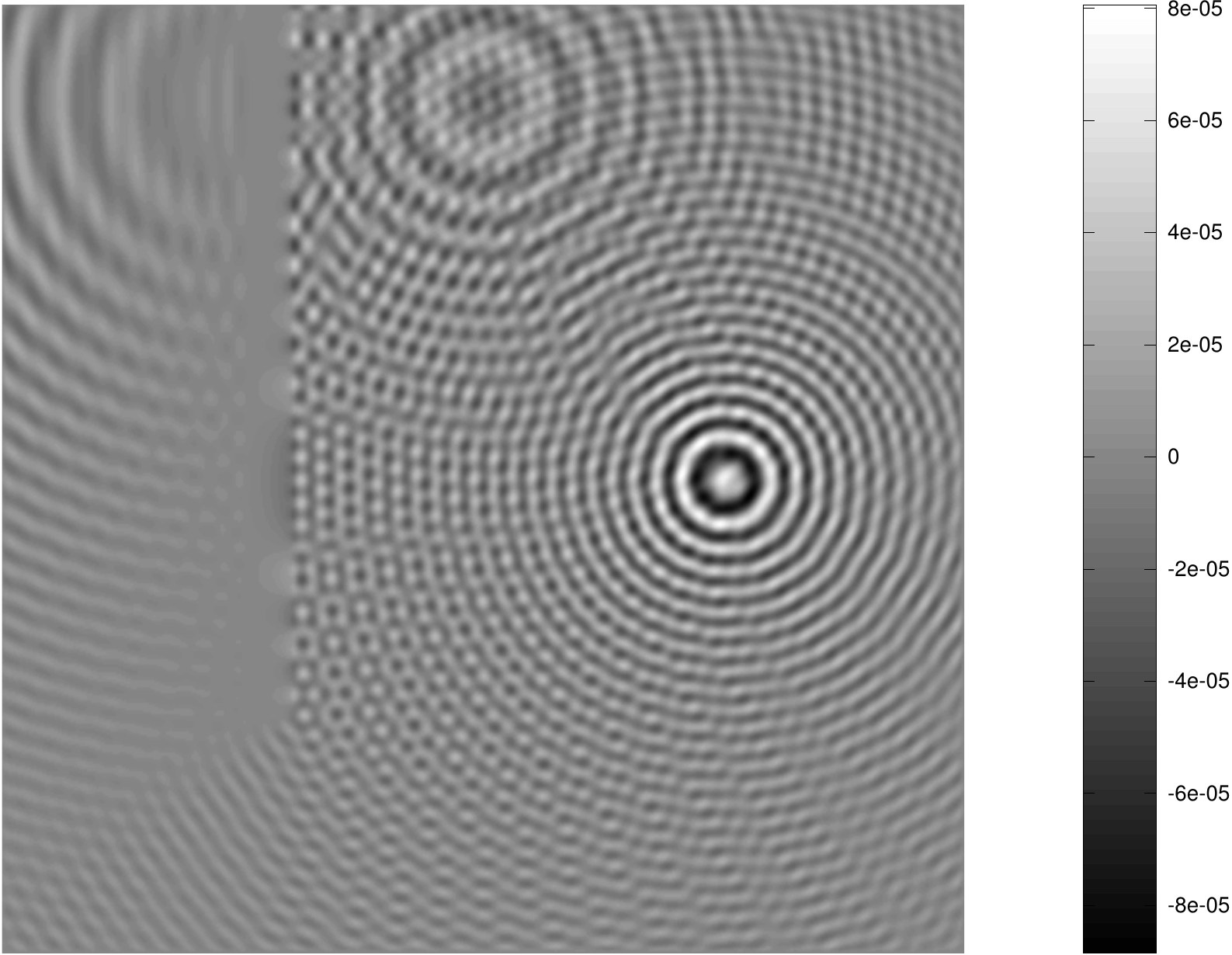} &
\includegraphics[width=2.5in]{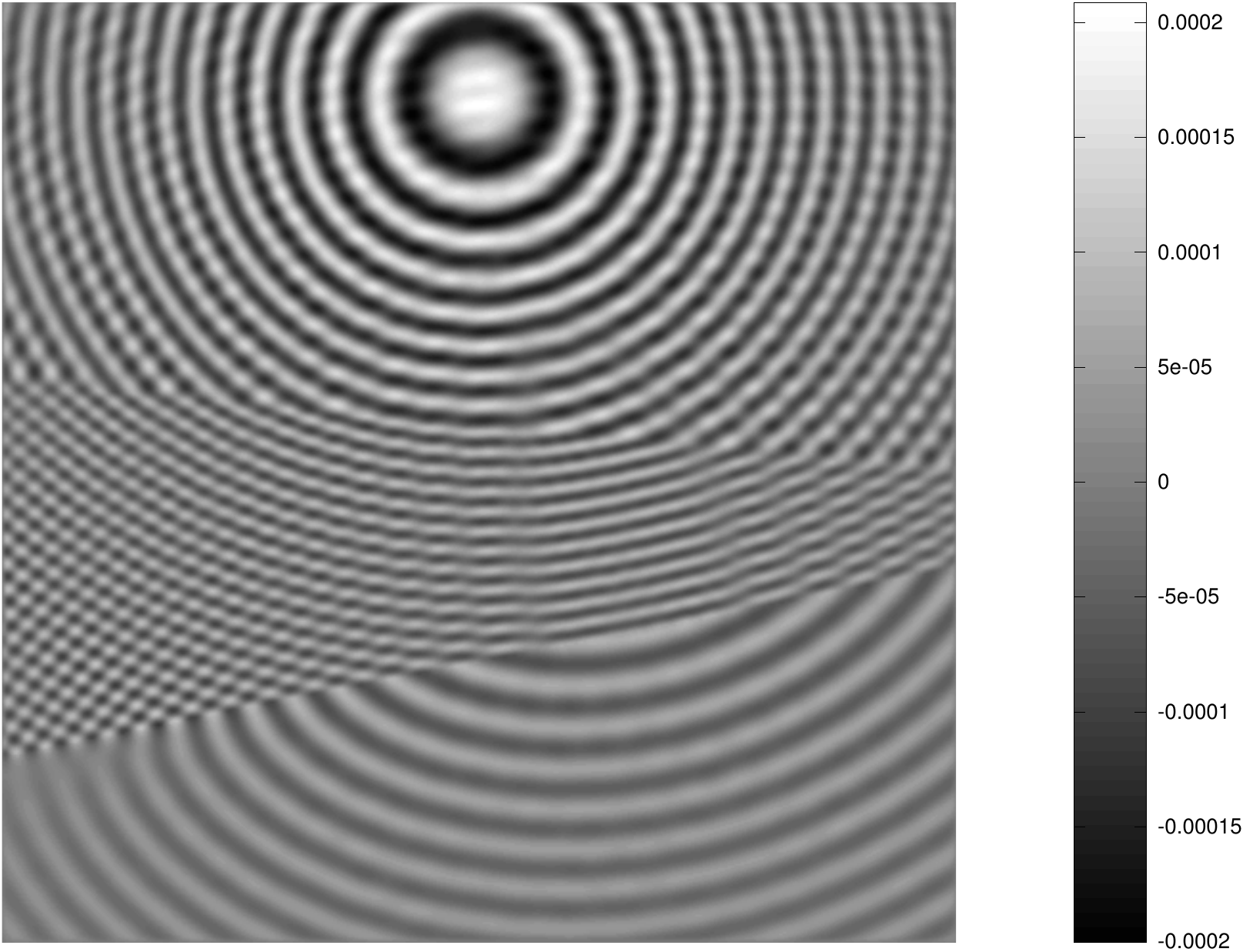} \\
\includegraphics[width=2.5in]{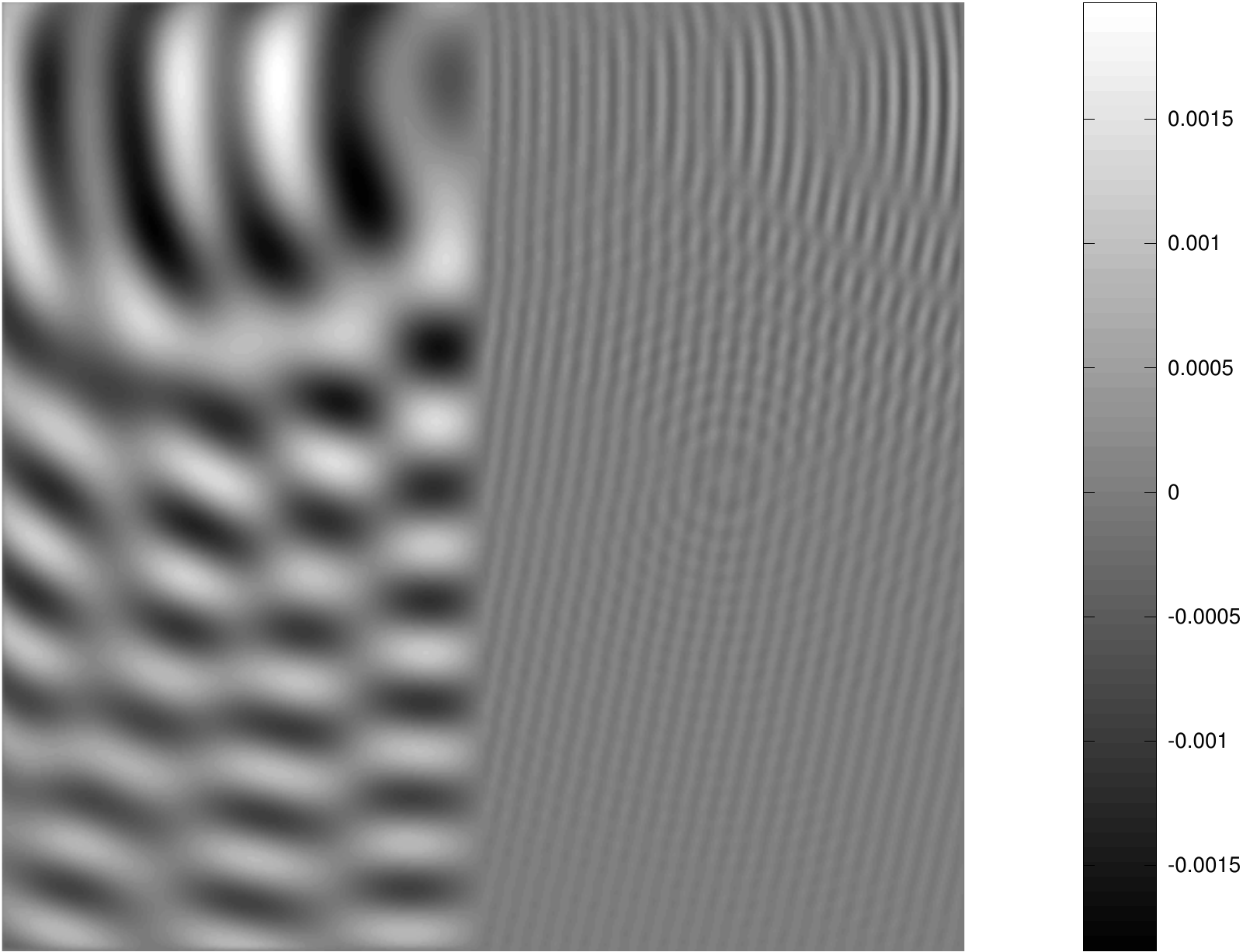} &
\includegraphics[width=2.5in]{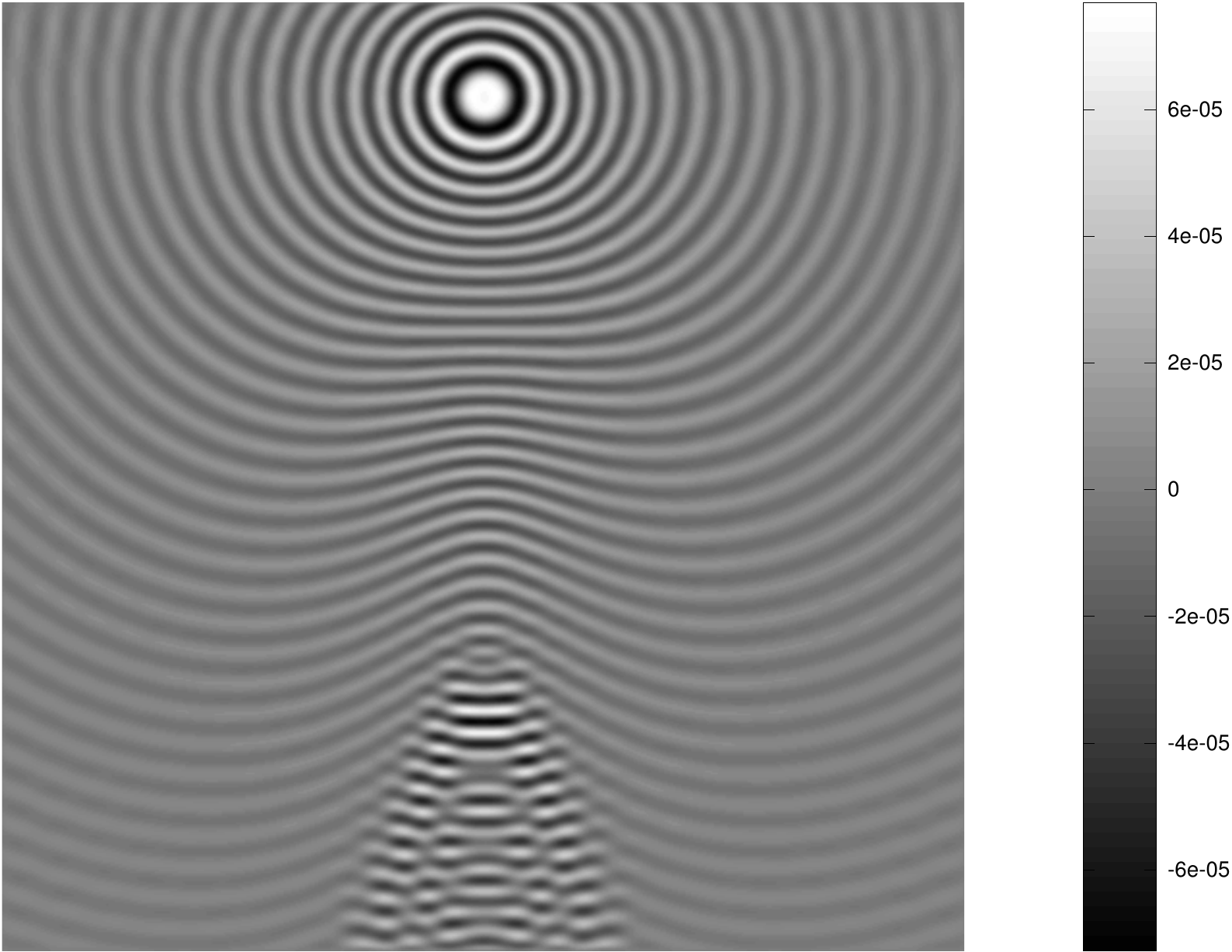}
\end{array}
$
\caption{A single $x_2 x_3$ plane from each of the four analytical velocity 
models over a $500^3$ grid with the smallest wavelength resolved with ten grid 
points.
(Top-left) the three-shot solution for the barrier model near $x_1=0.7$, 
(bottom-left) the three-shot solution for the two-layer model near $x_1=0.7$, 
(top-right) the single-shot solution for the wedge model near $x_1=0.7$, and
(bottom-right) the single-shot solution for the waveguide model near 
$x_1=0.55$.}
\end{figure}

It was also observed that, at least for the waveguide problem, 
the convergence rate would be significantly improved if 6 grid points of PML 
were used instead of 5. 
In particular, the 52 iterations reported in Table \ref{tbl:convergence} reduce 
to 27 if the PML size is increased by one. 
Interestingly, the same number of iterations
are required for convergence of the waveguide problem at half the frequency 
(and half the resolution) with five grid points of PML. Thus, it appears that 
the optimal PML size is a slowly growing function of the 
frequency.\footnote{A similar observation is also made 
in~\cite{Stolk-sweeping}.}
We also note that, though it was not intentional, each of the these first four 
velocity
models is invariant in one or more direction, and so it would be straightforward
to sweep in a direction such that only $O(1)$ panel factorizations would need
to be performed, effectively reducing the complexity of the setup phase to 
$O(\gamma^2 N)$.

The last experiment was meant to simultaneously test the convergence rates and 
scalability of the new sweeping preconditioner using the Overthrust velocity 
model (see Fig.~\ref{fig:overthrust}) at various frequencies, and the results
are reported in Table \ref{tbl:overthrust-test}.
It is important to notice that this is not a typical weak scaling test, as
the number of grid points per process was fixed, {\em not} the local memory 
usage or computational load, which both grow superlinearly with respect to the
total number of degrees of freedom. Nevertheless, including the setup phase, 
it took less than three minutes to solve the 3.175 Hz problem with 
four right-hand sides with 128 cores, and just under seven and a half minutes
to solve the corresponding 8 Hz problem using 2048 cores.
Also, while it is by no means the main 
message of this paper, the timings without making use of selective inversion 
are also reported in parentheses, as the technique is not widely 
implemented.\footnote{Other than Clique, the only other implementation appears 
to be in DSCPACK~\cite{Raghavan-dscpack}.}

\begin{table}
\footnotesize
\centering
\begin{tabular}{r|c|c|c|c|c|}
\cline{2-6}
 & \multicolumn{5}{|c|}{number of processes} \\
\cline{2-6}
 & 128 & 256 & 512 & 1024 & 2048  \\
\hline
\multicolumn{1}{|r|}{Hz}
 & 3.175
 & 4
 & 5.04
 & 6.35
 & 8 \\
\hline
\multicolumn{1}{|r|}{grid}
 & $319^2 \!\times\! 75$
 & $401^2 \!\times\! 94$
 & $505^2 \!\times\! 118$
 & $635^2 \!\times\! 145$
 & $801^2 \!\times\! 187$ \\
\hline
\multicolumn{1}{|r|}{setup (sec)}
 & 48.40 (46.23)
 & 66.33 (63.41) 
 & 92.95 (86.90)
 & 130.4 (120.6)
 & 193.0 (175.2) \\
\hline
\multicolumn{1}{|r|}{apply (sec/rhs)}
 & 0.468 (1.07)
 & 0.550 (1.28) 
 & 0.645 (2.40)
 & 0.700 (3.33)
 & 0.880 (6.13) \\
\hline
\multicolumn{1}{|r|}{3 digits (iter's)}
 & 42
 & 44
 & 42
 & 39
 & 40 \\
\hline
\multicolumn{1}{|r|}{4 digits (iter's)}
 & 54
 & 57
 & 57
 & 58
 & 58 \\
\hline
\multicolumn{1}{|r|}{5 digits (iter's)}
 & 63
 & 69
 & 70
 & 68
 & 72 \\
\hline
\end{tabular}
\caption{Convergence rates and timings on TACC's Lonestar for the 
SEG/EAGE Overthrust model, where timings in parentheses do not make use of 
selective inversion. All cases used a double-precision second-order 
stencil with five grid spacings for all PML (with an amplitude
of 7.5), and a damping parameter of $2.25 \pi$.
The preconditioner was configured with four planes per panel and eight 
processes per node, and the `apply' timings are with respect to a single 
application of the preconditioner to four right-hand sides.}
\label{tbl:overthrust-test}
\end{table}

\begin{figure}
\centering
$
\begin{array}{c}
\includegraphics[width=5in,height=1.1673in]{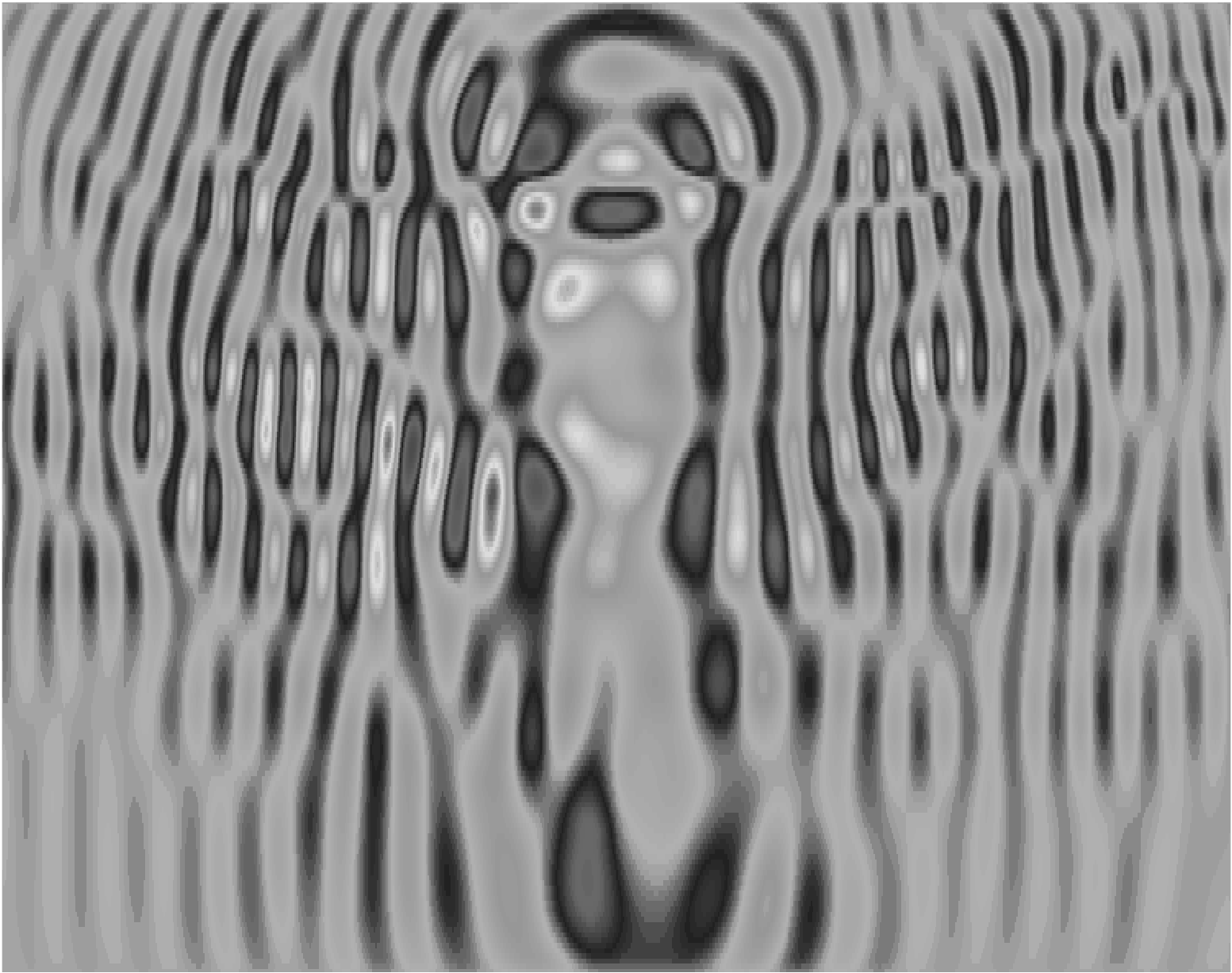}\\
\includegraphics[width=5in,height=1.1673in]{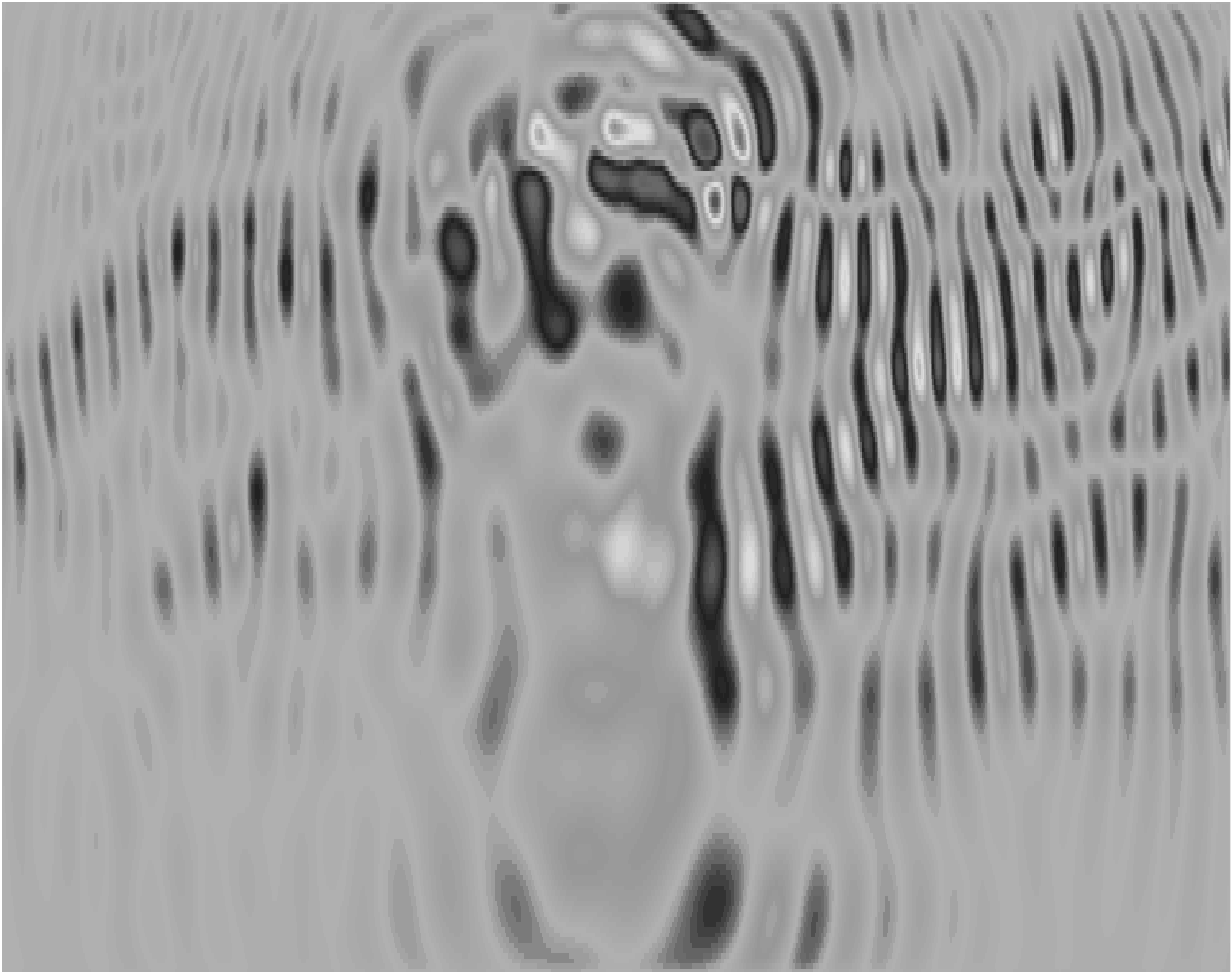}\\
\includegraphics[width=5in,height=5in]{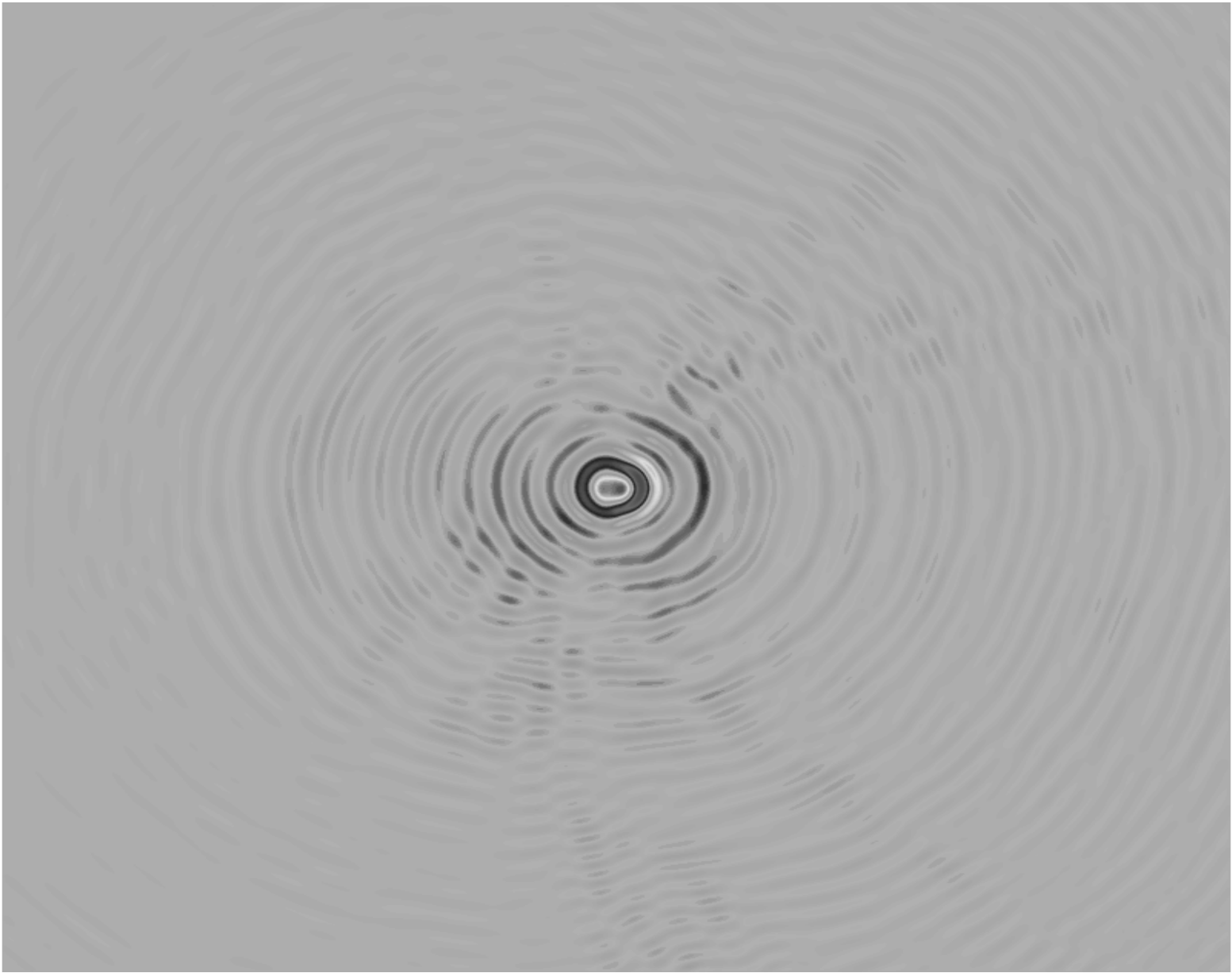}
\end{array}
$
\caption{Three planes from an 8 Hz solution with the Overthrust model at its 
native resolution, $801 \times 801 \times 187$, with a single localized shot 
at the center of the $x_1 x_2$ plane at a depth of 456 m: 
(top) a $x_2 x_3$ plane near $x_1=14$ km, 
(middle) an $x_1 x_3$ plane near $x_2=14$ km,
and (bottom) an $x_1 x_2$ plane near $x_3=0.9$ km.}
\label{fig:overthrust-planes}
\end{figure}

\section{Conclusions}
A parallelization of the {\em moving PML} sweeping preconditioner has been 
presented which has allowed us to efficiently solve 3D Helmholtz equations in 
parallel with essentially $O(1)$ iterations, with 
the only observed frequency-dependence arising from a moderate growth in the 
PML size with increasing frequency. This size of the PML, $\gamma(\omega)$ 
was explained to result in a linear growth in the memory requirements of the
preconditioner and a quadratic growth in the setup cost. Results were then 
presented for a variety of models, one of which had a velocity field which 
varied by ten orders of magnitude, and convergence was shown to be essentially
independent of frequency for the challenging Overthrust model.

Also, despite the requirement that each panel must be solved against one at a 
time when applying the preconditioner, a custom approach was introduced and 
implemented which eliminates most of the communication associated with 
performing traditional black-box sparse-direct factorizations and solves 
over each subdomain. These implementations are now released as part of the 
open-source packages Clique and Parallel Sweeping Preconditioner (PSP). There 
are at least five important directions for future work:
\begin{itemize}
\item developing a heuristic for tailoring the PML profile to the 
      velocity field,
\item extending the preconditioner to more general discretizations and time-harmonic wave equations,
\item finding a fast preconditioner for problems with large cavities
      (perhaps through more general local auxiliary problems), 
\item testing the performance improvements resulting from simultaneously 
      factoring the subdomain problems and then redistributing the 
      frontal trees, as well as a \verb!trsm! approach to solving many 
      right-hand sides, and
\item carefully studying the spectrum of the preconditioned operator for 
      various classes of velocity models.
\end{itemize}


%

\section*{Availability} The distributed dense linear algebra library, Elemental,
is available under the New BSD License at 
\url{http://code.google.com/p/elemental}. The distributed multifrontal solver, 
Clique, is available under the GPLv3 at 
\url{http://github.com/poulson/Clique}. The Parallel Sweeping 
Preconditioner (PSP) code is available under the GPLv3 at
\url{http://github.com/poulson/PSP}.

\section*{Acknowledgments}
The authors acknowledge TACC for usage of their computing resources and 
thank Bill Barth for suggesting the \verb!tacc_affinity! option and Tommy 
Minyard for helping with several large runs. We also thank Anshul Gupta for 
WSMP details and the referees for their insightful comments.

\end{document}